\documentclass[11pt]{article}
\usepackage[T1]{fontenc}
\usepackage{amsmath,amsthm,mathtools,mathdots}
\usepackage{newtxtext,newtxmath}
\usepackage{microtype}
\usepackage{enumitem}
\usepackage{geometry}
\usepackage{hyperref}
\usepackage{cleveref}
\usepackage{booktabs}
\hypersetup{colorlinks=true,linkcolor=black,citecolor=black,urlcolor=black}
\allowdisplaybreaks

\newtheorem{theorem}{Theorem}[section]
\newtheorem{proposition}[theorem]{Proposition}
\newtheorem{lemma}[theorem]{Lemma}
\newtheorem{corollary}[theorem]{Corollary}

\newtheorem{definition}[theorem]{Definition}

\crefname{theorem}{Theorem}{Theorems}
\Crefname{theorem}{Theorem}{Theorems}
\crefname{proposition}{Proposition}{Propositions}
\Crefname{proposition}{Proposition}{Propositions}
\crefname{lemma}{Lemma}{Lemmas}
\Crefname{lemma}{Lemma}{Lemmas}
\crefname{corollary}{Corollary}{Corollaries}
\Crefname{corollary}{Corollary}{Corollaries}
\crefname{remark}{Remark}{Remarks}
\Crefname{remark}{Remark}{Remarks}

\providecommand{\Bbbk}{\mathbb{k}}

\newcommand{\Z}{\mathbb Z}
\newcommand{\Aut}{\operatorname{Aut}}
\newcommand{\Hom}{\operatorname{Hom}}
\newcommand{\ev}{\mathrm{ev}}

\newcommand{\kk}{\Bbbk}

\title{Automorphisms Induced on Cartan Tori and Centers of\\
Quantized Enveloping Algebras}
\author{}
\date{}

\begin{document}

\begin{center}
\vspace*{1.5em}
{\LARGE\bfseries Automorphisms Induced on Cartan Tori and Centers of\\[0.3em]
Quantized Enveloping Algebras\par}
\end{center}
\vspace{1.5em}

\begin{center}
	Yu Wang\\[4pt]
	{\small
		School of Mathematics and Physics, Jiangsu University of Technology,\\
		Changzhou, Jiangsu 213001, China\\
		E-mail: \texttt{yu.wang@mail.sdu.edu.cn}\par}
\end{center}

\begin{abstract}

Let $U=U_q(\mathfrak g)$ be the root-lattice form of the quantized
enveloping algebra associated with a complex simple Lie algebra
$\mathfrak g$, over a field $\Bbbk$ of characteristic zero.
The field $\Bbbk$ need not be algebraically closed, and we assume that
$q\in\Bbbk^\times$ is not a root of unity.
We determine the groups of automorphisms induced by algebra
automorphisms of $U$ on its Cartan torus $U^0$ and its center $Z(U)$.
These groups are
\[
\begin{aligned}
	\Aut_U(U^0)&\cong\Hom(Q,\{\pm1\})\rtimes\Aut(\Phi),\\
	\Aut_U(Z(U))&\cong\Hom((P/Q)[2],\{\pm1\})\rtimes\Gamma,
\end{aligned}
\]
where $Q$ and $P$ are the root and weight lattices, $\Phi$ is the root
system, $\Gamma$ is the Dynkin-diagram automorphism group of $\mathfrak g$,
and $(P/Q)[2]$ is the $2$-torsion subgroup of $P/Q$.
We give explicit descriptions of these groups for every simple Lie
type, including the exceptional types.
Writing $\rho_0$ and $\rho_Z$ for the restriction maps to $U^0$ and
$Z(U)$, respectively, we prove that $\ker\rho_0\subseteq\ker\rho_Z$
and determine the quotient $\ker\rho_Z/\ker\rho_0$.
The proof uses characters of finite-dimensional modules and Weyl-group
invariant fields to show that the lattice automorphism determined by
the toral restriction preserves the root system.
Central characters and the quantum Harish--Chandra isomorphism then
determine the induced action on the center without choosing explicit
central generators.
\end{abstract}

\section{Introduction}

Let $U=U_q(\mathfrak g)$ be the root-lattice form of the quantized enveloping algebra of a complex simple Lie algebra $\mathfrak g$. We assume that the ground field $\Bbbk$ has characteristic zero but is not necessarily algebraically closed, and that $q\in\Bbbk^\times$ is not a root of unity. All automorphisms are understood to be $\Bbbk$-algebra automorphisms. Standard references for the structure and representation theory of $U$ are [6, 9, 15]. Determining $\Aut(U)$ is difficult in general: arbitrary algebra automorphisms need not preserve the Hopf algebra structure, and even a prescribed action on the Cartan torus $U^0$ does not determine the images of the generators $E_i$ and $F_i$. We study the restrictions of these automorphisms to 
$U^0$ and $Z(U)$, and determine the images of the two
restriction maps in terms of the root system and its lattices.

Twietmeyer [18] described the Hopf automorphisms of $U$ in terms of diagonal and Dynkin-diagram automorphisms. Chin and Musson [7] obtained a description of the Hopf automorphism group over $\mathbb Q(q)$, with $q$ an indeterminate. Lusztig's braid-group operators [15] are algebra automorphisms whose restrictions to $U^0$ realize the Weyl-group action. In rank one, Alev and Chamarie [1] determined the full algebra automorphism group of $U_q(\mathfrak{sl}_2)$. Alev and Dumas [2, \S3.1] recall this classification in a weight-lattice presentation, which differs from the root-lattice form considered here. For the positive part $U_q^+(\mathfrak g)$ in rank at least two, Andruskiewitsch and Dumas [3] conjectured that every algebra automorphism is a composition of rescalings of the generators and a Dynkin-diagram automorphism. Launois [12] determined this group in type $B_2$. Yakimov [19] proved the conjecture using rigidity of completed quantum tori, obtaining $\Aut(U_q^+(\mathfrak g))\cong(\Bbbk^\times)^n\rtimes\Gamma$, where $n$ is the rank of $\mathfrak g$ and $\Gamma$ is the
automorphism group of its Dynkin diagram. Since automorphisms of $U$ need not preserve $U_q^+(\mathfrak g)$, this result does not directly determine their restrictions to $U^0$ or $Z(U)$.

With $q$ an indeterminate, Zha [20] studied deformed automorphisms of $U_q(\mathfrak g)$, which preserve an integral form and induce Lie algebra automorphisms of $\mathfrak g$ in the classical limit $q=1$. He proved that their induced automorphisms of the root lattice preserve the root system. We establish this constraint for every algebra automorphism under the assumptions above, without imposing deformation conditions.

The description of the units of $U$ implies that every automorphism
preserves $U^0$. Since every automorphism also preserves $Z(U)$,
restriction defines group homomorphisms
\[
\begin{aligned}
	\rho_0&:\Aut(U)\longrightarrow\Aut(U^0),
	&\sigma&\longmapsto\sigma|_{U^0},\\
	\rho_Z&:\Aut(U)\longrightarrow\Aut(Z(U)),
	&\sigma&\longmapsto\sigma|_{Z(U)}.
\end{aligned}
\]
We denote their images by $\Aut_U(U^0)$ and $\Aut_U(Z(U))$,
respectively.

Let $\Phi$ be the root system, with root lattice $Q$, weight lattice $P$, and Weyl group $W$. Theorem~4.1 shows that every $\sigma\in\Aut(U)$ satisfies
\[
\sigma(K_\lambda)=\varepsilon(\lambda)K_{\varphi(\lambda)}
\qquad(\lambda\in Q),
\]
where $\varepsilon\in\Hom(Q,\{\pm1\})$ and $\varphi\in\Aut(\Phi)=W\rtimes\Gamma$. Conversely, every such toral action is realized by composing sign automorphisms, Lusztig's braid-group operators, and graph automorphisms. Hence
\[
\Aut_U(U^0)\cong\Hom(Q,\{\pm1\})\rtimes\Aut(\Phi).
\]
Theorem~4.10 shows that an automorphism of $Z(U)$ extends to
an algebra automorphism of $U$ if and only if it is induced
by a composition of sign and graph automorphisms.
Consequently,
\[
\Aut_U(Z(U))\cong\Hom((P/Q)[2],\{\pm1\})\rtimes\Gamma,
\]
where $(P/Q)[2]$ is the $2$-torsion subgroup of $P/Q$. The character-group factor is obtained by restricting sign characters
to $Q\cap2P$ and using the natural isomorphism $(Q\cap2P)/2Q\cong(P/Q)[2]$.

The two restrictions are related by
$\ker\rho_0\subseteq\ker\rho_Z$.
Thus automorphisms that agree on $U^0$ also agree on $Z(U)$, and $\rho_Z$ factors through $\rho_0$. Proposition~4.12 characterizes the toral restrictions of
automorphisms fixing the center pointwise and gives
\[
\ker\rho_Z/\ker\rho_0
\cong\Hom\bigl(Q/(Q\cap2P),\{\pm1\}\bigr)\rtimes W,
\]
where $W$ acts on characters by $(w\cdot\delta)(\lambda)=\delta(w^{-1}\lambda)$. The determination of $\ker\rho_0$ and of the full automorphism group
remains outside the scope of this paper.

The proofs use finite-dimensional irreducible type-$1$ modules. After removing the sign part of an automorphism, we study its effect on weights and characters through twisting. Together with Weyl-group invariant fields, this yields the root-system constraint. We then compare central characters and apply the quantum Harish--Chandra isomorphism to determine the induced action on the center, without choosing explicit central generators for each Cartan type. For the structure of the center, we use [10, 11, 16, 17]; see also [4, 13, 14].

The paper is organized as follows. Section~2 collects the preliminaries. Section~3 establishes the root-system constraint, and Section~4 determines the two restriction images and relates their kernels. Root-system numbering follows Bourbaki [5]; the Cartan-matrix convention is specified in Section~2.

\section{Preliminary results}
Throughout the paper, let $\Bbbk$ be a field of characteristic $0$, and let
$\mathfrak g$ be a finite-dimensional simple Lie algebra over $\mathbb C$ of
rank $n\ge1$, with root system $\Phi$, simple system
$\Delta=\{\alpha_1,\ldots,\alpha_n\}$, and Cartan matrix
\[
A=(a_{ij})_{1\le i,j\le n}.
\]
Fix the standard Weyl-group-invariant symmetric bilinear form
$(\, ,\,)$ on the rational span of $\Phi$.
Normalize this form so that the numbers
\[
\frac{(\alpha_i,\alpha_i)}{2}\qquad(1\leq i\leq n)
\]
are positive integers with greatest common divisor $1$. We use Jantzen's convention~\cite{Jantzen}
\[
a_{ij}=\langle\alpha_j,\alpha_i^\vee\rangle
=\frac{2(\alpha_i,\alpha_j)}{(\alpha_i,\alpha_i)}.
\]
In the nonsimply-laced cases, this is the transpose of the Cartan-matrix convention used in \cite{Humphreys}; all Cartan matrices below follow the convention displayed above. We put
\[
d_i=\frac{(\alpha_i,\alpha_i)}{2},\qquad q_i=q^{d_i},
\qquad
D=\operatorname{diag}\left(
d_1,\ldots, d_n
\right).
\]
Then $DA$ is symmetric. We assume throughout that $q\in\Bbbk^\times$ is not a root of unity. 
\subsection{The quantized enveloping algebra and its center}
Let $U=U_q(\mathfrak g)$ be the root-lattice form of the Drinfeld--Jimbo quantized enveloping algebra, generated by $E_i,F_i,K_i^{\pm1}$ $(1\le i\le n)$, with the usual defining relations
\[
\begin{array}{@{}c@{\qquad\qquad}c@{}}
K_iK_j=K_jK_i,
&
K_iK_i^{-1}=K_i^{-1}K_i=1,
\\[2pt]
K_iE_jK_i^{-1}=q_i^{a_{ij}}E_j,
&
K_iF_jK_i^{-1}=q_i^{-a_{ij}}F_j,
\\[5pt]
\multicolumn{2}{c}{
E_iF_j-F_jE_i
=\delta_{ij}\dfrac{K_i-K_i^{-1}}{q_i-q_i^{-1}}.
}
\end{array}
\]
For $i\ne j$, the quantum Serre relations are
\[
\sum_{k=0}^{1-a_{ij}}(-1)^k
\begin{bmatrix}1-a_{ij}\\ k\end{bmatrix}_{q_i}
E_i^{\,1-a_{ij}-k}E_jE_i^k=0,
\qquad
\sum_{k=0}^{1-a_{ij}}(-1)^k
\begin{bmatrix}1-a_{ij}\\ k\end{bmatrix}_{q_i}
F_i^{\,1-a_{ij}-k}F_jF_i^k=0,
\]
where
\[
[m]_{q_i}=\frac{q_i^m-q_i^{-m}}{q_i-q_i^{-1}},
\qquad
[m]_{q_i}!=\prod_{s=1}^{m}[s]_{q_i},
\qquad
\begin{bmatrix}m\\ k\end{bmatrix}_{q_i}
=\frac{[m]_{q_i}!}{[k]_{q_i}![m-k]_{q_i}!}.
\]

Let 
\[Q=\bigoplus_{i=1}^{n}\Z\alpha_i, \qquad P=\bigoplus_{i=1}^{n}\Z\omega_i\]
be the root and weight lattices, where $\langle\omega_i,\alpha_j^\vee\rangle=\delta_{ij}$. Let
\[
Q^+=\sum_{i=1}^n\mathbb Z_{\geq0}\alpha_i,
\qquad
P^+=\sum_{i=1}^n\mathbb Z_{\geq0}\omega_i
\]
be the positive root cone and the set of dominant integral
weights, respectively. We use the dominance order on $P$:
$\mu\leq\lambda$ if $\lambda-\mu\in Q^+$.
We write $\mu<\lambda$ when $\mu\leq\lambda$ and $\mu\neq\lambda$.
For $\lambda=\sum_i m_i\alpha_i\in Q$, write $K_\lambda=K_1^{m_1}\cdots K_n^{m_n}$. Then
\[
U^0=\Bbbk[K_1^{\pm1},\ldots,K_n^{\pm1}]=\Bbbk[Q].
\]
Let $\mathcal E:=\Hom(Q,\{\pm1\})$ be the group of
homomorphisms from the additive group $Q$ to the multiplicative
group $\{\pm1\}\subseteq\Bbbk^\times$, with pointwise multiplication.
Let $\Gamma:=\Aut(\operatorname{Dyn}(\mathfrak g))$ be the
Dynkin-diagram automorphism group of $\mathfrak g$. Since $Q$ is free of rank $n$,
\[
\mathcal{E}\cong(\Z/2\Z)^n,
\qquad |\mathcal{E}|=2^n.
\]
Let $U^+$ and $U^-$ be the subalgebras generated by the $E_i$
and the $F_i$, respectively. We use the $Q$-grading given by $\deg K_i^{\pm1}=0$, $\deg E_i=\alpha_i$, and $\deg F_i=-\alpha_i$, together with the triangular decomposition $U=U^-U^0U^+$. For $\nu\in Q^+$, let $U_\nu^+$ and $U_{-\nu}^-$ denote
the homogeneous components of $U^+$ and $U^-$ of degrees
$\nu$ and $-\nu$, respectively.

Let $U_0$ denote the zero component for the root-lattice grading, to be distinguished from the toral subalgebra $U^0$. The triangular decomposition gives
\[
U_0
=
U^0 \oplus
\bigoplus_{\nu\in Q^{+}\setminus\{0\}}
U^-_{-\nu}U^0U^+_\nu.
\]
Let $\pi:U_0\to U^0$ be projection onto the first summand. Since $q$ is not a root of unity, $Z(U)\subseteq U_0$. Let $\rho$ be half the sum of the positive roots, and define
\[
\gamma_{-\rho}(K_\lambda)=q^{-(\rho,\lambda)}K_\lambda,
\qquad
U^0_{\ev}=\bigoplus_{\lambda\in Q\cap2P}\kk K_\lambda.
\]
Since $(\rho,\alpha_i)=d_i\in\mathbb Z$, we have
$(\rho,\lambda)\in\mathbb Z$ for every $\lambda\in Q$.
Thus all the scalars in this definition belong to $\Bbbk$. The action of $W$ on $U^0$ is the action $w(K_\lambda)=K_{w\lambda}$.

\begin{proposition}\label{2.1}
The map
\[
\gamma_{-\rho}\circ\pi \colon
Z(U)\xrightarrow{\sim}(U_{\mathrm{ev}}^0)^W
\]
is an isomorphism of \(\Bbbk\)-algebras.
\end{proposition}
\begin{proof}
This is the quantum Harish--Chandra isomorphism;
see \cite[Theorem 6.25 and \S6.26]{Jantzen},
together with \S8.30 for the extension to arbitrary parameters that are not roots of unity.
\end{proof}

\subsection{Units, toral actions, and exponent matrices}
We first recall the standard description of units.

\begin{lemma}\label{lem:units}
The units of \(U\) are precisely
\[
U^\times
=
\{cK_\lambda :
c\in\Bbbk^\times,\ \lambda\in Q\}.
\]
\end{lemma}

\begin{proof}
The argument of \cite[Lemma~2]{Zha} applies here.
It uses the root-lattice grading, the trivial module,
and the leading-term multiplication of PBW monomials,
all of which remain valid under our assumptions.
In particular, it does not require \(q\) to be transcendental.
\end{proof}

\begin{lemma}\label{lem:monomial}
	Let $\sigma\in\Aut(U)$. Then there exist unique elements
	$\varepsilon\in\mathcal E$ and $\varphi\in\Aut_{\mathbb Z}(Q)$ such that
	\[
	\sigma(K_\lambda)=\varepsilon(\lambda)K_{\varphi(\lambda)}
	\qquad(\lambda\in Q).
	\]
\end{lemma}
\begin{proof}
	By Lemma~2.2, every unit of $U$ has a unique expression $aK_\lambda$,
	where $a\in\Bbbk^\times$ and $\lambda\in Q$. Since $\sigma$ maps units
	to units, there are uniquely determined maps
	$c:Q\to\Bbbk^\times$ and $\varphi:Q\to Q$ such that
	\[
	\sigma(K_\lambda)=c(\lambda)K_{\varphi(\lambda)}
	\qquad(\lambda\in Q).
	\]
	For $\lambda,\mu\in Q$, applying $\sigma$ to
	$K_{\lambda+\mu}=K_\lambda K_\mu$ gives
	\[
	c(\lambda+\mu)K_{\varphi(\lambda+\mu)}
	=c(\lambda)c(\mu)K_{\varphi(\lambda)+\varphi(\mu)}.
	\]
	Uniqueness of the expression for a unit therefore implies
	\[
	\varphi(\lambda+\mu)=\varphi(\lambda)+\varphi(\mu),
	\qquad c(\lambda+\mu)=c(\lambda)c(\mu).
	\]
	Thus $\varphi$ is additive and $c$ is a character of $Q$.
	
	Applying the same argument to $\sigma^{-1}$, write
	\[
	\sigma^{-1}(K_\nu)=d(\nu)K_{\psi(\nu)}
	\qquad(\nu\in Q),
	\]
	where $d:Q\to\Bbbk^\times$ is a character and $\psi:Q\to Q$ is additive.
	Since $\sigma^{-1}$ fixes scalars, we have
	\[
	K_\lambda=\sigma^{-1}\bigl(\sigma(K_\lambda)\bigr)
	=c(\lambda)d\bigl(\varphi(\lambda)\bigr)
	K_{\psi(\varphi(\lambda))}.
	\]
	Comparing exponents yields $\psi(\varphi(\lambda))=\lambda$.
	Similarly, $\sigma(\sigma^{-1}(K_\nu))=K_\nu$ gives
	$\varphi(\psi(\nu))=\nu$. Hence $\psi=\varphi^{-1}$ and
	$\varphi\in\Aut_{\mathbb Z}(Q)$.
	
	Finally, let $\eta:U\to\Bbbk$ be the algebra homomorphism defined by
	\[
	\eta(E_i)=\eta(F_i)=0,\qquad \eta(K_i)=1.
	\]
	Set $\chi=\eta\circ\sigma$. Then
	$\chi(K_i)=c(\alpha_i)$ and $\chi(K_i^{-1})=c(\alpha_i)^{-1}$.
	Since $\Bbbk$ is commutative, $\chi(E_iF_i-F_iE_i)=0$.
	Applying $\chi$ to
	\[
	E_iF_i-F_iE_i=\frac{K_i-K_i^{-1}}{q_i-q_i^{-1}}
	\]
	gives
	\[
	0=\frac{c(\alpha_i)-c(\alpha_i)^{-1}}{q_i-q_i^{-1}}.
	\]
	Since $q$ is not a root of unity, $q_i-q_i^{-1}\ne0$.
	Thus $c(\alpha_i)^2=1$, and hence $c(\alpha_i)\in\{\pm1\}$ for every $i$.
	Since the simple roots form a
	$\mathbb Z$-basis of $Q$, the character $c$ takes values in $\{\pm1\}$.
	Setting $\varepsilon=c$ proves the assertion, including uniqueness.
\end{proof}

Since every algebra automorphism of $U$ preserves $U^0$ by Lemma~\ref{lem:monomial}, we may write
\[
\Aut_{U}(U^0)
:=\{\sigma|_{U^0}:\sigma\in\Aut(U)\}.
\]
By Lemma~\ref{lem:monomial}, every $\sigma\in\Aut(U)$ acts on the toral generators by signed Laurent monomials. Writing
\begin{equation*}\label{eq:sigma-Ki}
\sigma(K_i)=\varepsilon_iK_1^{b_{i1}}\cdots K_n^{b_{in}},
\qquad \varepsilon_i\in\{\pm1\},\qquad 1\le i\le n,
\end{equation*}
we define the associated \emph{exponent matrix}
$A_\sigma:=(b_{ij})_{1\le i,j\le n}.$ Since \(\varphi\in\operatorname{Aut}_{\mathbb Z}(Q)\) and
\(\varphi(\alpha_i)=\sum_{j=1}^n b_{ij}\alpha_j\),
the transpose \(A_\sigma^{T}\) is the matrix of
\(\varphi\) with respect to the basis
\(\alpha_1,\ldots,\alpha_n\).
Hence \(A_\sigma\in \mathrm{GL}_n(\mathbb Z)\).

For $\delta\in\mathcal{E}$, put
$\delta_i=\delta(\alpha_i)$ and define the sign automorphism
$s_\delta$ of $U$ by
\[
s_\delta(E_i)=\delta_iE_i,\qquad
s_\delta(F_i)=F_i,\qquad
s_\delta(K_i)=\delta_iK_i
\quad(1\leq i\leq n).
\]
In particular, $s_\delta(K_\lambda)=\delta(\lambda)K_\lambda$
for every $\lambda\in Q$.

We call an automorphism $\tau$ \emph{sign-free on the Cartan torus}
if $\tau(K_\lambda)=K_{\varphi(\lambda)}$ for every $\lambda\in Q$,
where $\varphi\in\operatorname{Aut}_{\mathbb Z}(Q)$. 
\begin{lemma}\label{lem:remove-sign}
	Let \(\sigma\in\operatorname{Aut}(U)\), and write
	\[
	\sigma(K_\lambda)=\varepsilon(\lambda)K_{\varphi(\lambda)}
	\qquad (\lambda\in Q)
	\]
	as in Lemma~2.3. Put
	\[
	\delta=\varepsilon\circ\varphi^{-1}\in\mathcal E,
	\qquad
	\tau=s_\delta\circ\sigma.
	\]
	Then $\tau$ is sign-free.
	In particular, \(A_\tau=A_\sigma\).
\end{lemma}

\begin{proof}
	Since \(\varphi\) is a lattice automorphism,
	\(\delta=\varepsilon\circ\varphi^{-1}\) belongs to
	\(\mathcal E\), and
	\(\delta(\varphi(\lambda))=\varepsilon(\lambda)\)
	for every \(\lambda\in Q\). Hence
$$\tau(K_\lambda)=s_\delta\bigl(\varepsilon(\lambda)K_{\varphi(\lambda)}\bigr)=\varepsilon(\lambda)\delta(\varphi(\lambda))
		K_{\varphi(\lambda)}
		=\varepsilon(\lambda)^2K_{\varphi(\lambda)}
		=K_{\varphi(\lambda)}.$$
Thus \(\tau\) and \(\sigma\) have the same exponent matrix.
\end{proof}

For $\lambda\in P$, let $[\lambda]$ denote its coordinate column with
respect to the ordered basis $\omega_1,\ldots,\omega_n$.
The natural action of $W$ on $P$ defines a faithful representation
\[
\vartheta\colon W\longrightarrow \mathrm{GL}_n(\mathbb Z),
\qquad
\vartheta(w)=\bigl([w(\omega_1)],\ldots,[w(\omega_n)]\bigr).
\]
Thus
\[
[w(\lambda)]=\vartheta(w)[\lambda]
\qquad (w\in W,\ \lambda\in P),
\]
and we write $\vartheta(W)=\{\vartheta(w):w\in W\}$.

Each $\xi\in\Gamma$ acts linearly on the root and coroot spaces by
\[
\xi(\alpha_i)=\alpha_{\xi(i)},
\qquad
\xi(\alpha_i^\vee)=\alpha_{\xi(i)}^\vee
\quad (1\le i\le n).
\]
Since $\xi$ preserves the Cartan matrix, these actions preserve the
natural pairing. Consequently,
\[
\langle\xi(\omega_i),\alpha_{\xi(j)}^\vee\rangle
=\langle\omega_i,\alpha_j^\vee\rangle
=\delta_{ij},
\]
so $\xi(\omega_i)=\omega_{\xi(i)}$. In particular, $\xi$ preserves $P$.
Let $[\xi]\in \mathrm{GL}_n(\mathbb Z)$ denote the corresponding permutation
matrix with respect to the ordered basis $\omega_1,\ldots,\omega_n$.
Then
\begin{equation*}\label{eq:xi-weight-action}
	[\xi(\lambda)]=[\xi][\lambda]
	\qquad (\lambda\in P).
\end{equation*}

\subsection{Standard automorphisms and root-system automorphisms}
We recall the standard automorphisms of $U$ and the structure
of $\Aut(\Phi)$. We then describe the subgroup of $\Aut(U^0)$
induced by sign automorphisms, Lusztig's braid-group operators,
and graph automorphisms.
\begin{lemma}\label{lem:standard-aut}
	The algebra $U=U_q(\mathfrak g)$ admits the following automorphisms.
	\begin{enumerate}
		\renewcommand{\labelenumi}{(\alph{enumi})}
		\item The involutive automorphism $I$ defined by
		\[
		I(E_i)=F_i,\qquad I(F_i)=E_i,\qquad I(K_i)=K_i^{-1}.
		\]
		\item For $a=(a_1,\ldots,a_n)\in(\Bbbk^\times)^n$, the diagonal
		automorphism $D_a$ defined by
		\[
		D_a(E_i)=a_iE_i,\qquad D_a(F_i)=a_i^{-1}F_i,\qquad D_a(K_i)=K_i.
		\]
		\item For $\xi\in\Gamma$, the graph automorphism $G_\xi$  defined by
		\[
		G_\xi(E_i)=E_{\xi(i)},\qquad G_\xi(F_i)=F_{\xi(i)},\qquad
		G_\xi(K_i)=K_{\xi(i)}.
		\]
		\item The sign automorphisms $s_\delta$, $\delta\in\mathcal E$, defined above.
		They satisfy
		\[
		s_\delta\circ s_\varepsilon=s_{\delta\varepsilon},\qquad
		s_\delta^2=\operatorname{id}_U
		\quad(\delta,\varepsilon\in\mathcal E),
		\]
		where $(\delta\varepsilon)(\lambda)=\delta(\lambda)\varepsilon(\lambda)$.
		The map $\delta\mapsto s_\delta$ identifies $\mathcal E$ with the
		subgroup of sign automorphisms of $U$.
		\item The Lusztig braid automorphisms $T_1,\ldots,T_n$, which satisfy
		\[
		T_i(K_\lambda)=K_{s_i(\lambda)}
		\qquad (\lambda\in Q),
		\]
		where $s_i$ is the simple reflection corresponding to $\alpha_i, 1\le i\le n$.
	\end{enumerate}
\end{lemma}

\begin{proof}
	For (a)--(d), the assignments preserve the defining relations of $U$,
	and their inverses are $I$, $D_{a^{-1}}$, $G_{\xi^{-1}}$, and $s_\delta$,
	respectively, where $a^{-1}=(a_1^{-1},\ldots,a_n^{-1})$.
	The identities in (d) follow by evaluation on the generators.
	Moreover, $s_\delta=\operatorname{id}_U$ implies
	$\delta(\alpha_i)=1$ for every $i$, and hence $\delta=1$.
	For (e), see \cite[Chapter~8]{Jantzen}.
\end{proof}
For later reference, the Dynkin-diagram automorphism groups for the simple types are listed below. We use the Bourbaki numbering of the simple roots for types $D$ and $E$.
\[
\begin{array}{c|c|l}
\text{type} & \Gamma=\Aut(\mathrm{Dyn}(\mathfrak g)) & \text{nontrivial Dynkin-diagram action}\\
\hline
A_1 & 1 & \text{none}\\
A_n\ (n\ge2) & \Z/2\Z & \xi(i)=n+1-i\\
B_n,\ C_n\ (n\ge2) & 1 & \text{none}\\
D_4 & S_3 & \text{arbitrary permutations of }\{1,3,4\}\text{, with }2\text{ fixed}\\
D_n\ (n\ge5) & \Z/2\Z & \xi(n-1)=n,\quad \xi(n)=n-1\\
E_6 & \Z/2\Z & \xi=(1\ 6)(3\ 5),\quad \xi(2)=2,\ \xi(4)=4\\
E_7,\ E_8,\ F_4,\ G_2 & 1 & \text{none}
\end{array}
\]

We identify $\Gamma$ with the subgroup of $\Aut(\Phi)$ preserving
the chosen simple system $\Delta$, via its action on the simple
roots. We recall the following standard result; see~[8, \S12.2].

\begin{lemma}\label{prop:root-aut-general}
	The root-system automorphism group is the semidirect product
	\[
	\Aut(\Phi)=W\rtimes\Gamma.
	\]
	In particular, every $\psi\in\Aut(\Phi)$ has a unique expression
	$\psi=w\xi$ with $w\in W$ and $\xi\in\Gamma$.
\end{lemma}

For $\varepsilon\in\mathcal E$ and $\varphi\in\Aut(\Phi)$,
let $\tau_{\varepsilon,\varphi}\in\Aut(U^0)$ be defined by
\[
\tau_{\varepsilon,\varphi}(K_\lambda)
=\varepsilon(\lambda)K_{\varphi(\lambda)}
\qquad(\lambda\in Q).
\]

\begin{proposition}\label{prop:cartan-realization-general}
	The automorphisms $\tau_{\varepsilon,\varphi}$, where
	$\varepsilon\in\mathcal E$ and $\varphi\in\Aut(\Phi)$,
	form a subgroup of $\Aut_U(U^0)$ isomorphic to
	$\mathcal E\rtimes\Aut(\Phi)$, with action
	\[
	\varphi\cdot\varepsilon=\varepsilon\circ\varphi^{-1}.
	\]
\end{proposition}

\begin{proof}
	By Lemma~2.6, write $\varphi=w\xi$, where $w\in W$ and $\xi\in\Gamma$.
	Choose a reduced expression $w=s_{i_1}\cdots s_{i_r}$ and set
	$T_w=T_{i_1}\circ\cdots\circ T_{i_r}$.
	By Lemma~2.5, $T_w(K_\lambda)=K_{w(\lambda)}$, and hence
	\[
	(T_w\circ G_\xi\circ s_\varepsilon)(K_\lambda)
	=\varepsilon(\lambda)K_{w\xi(\lambda)}
	=\tau_{\varepsilon,\varphi}(K_\lambda)
	\qquad(\lambda\in Q).
	\]
	Thus every $\tau_{\varepsilon,\varphi}$ extends to an algebra
	automorphism of $U$.
	
	The group operation is composition. For $\varepsilon,\delta\in\mathcal E$
	and $\varphi,\psi\in\Aut(\Phi)$, direct evaluation on $K_\lambda$ gives
	\[
	\tau_{\varepsilon,\varphi}\circ\tau_{\delta,\psi}
	=\tau_{\delta(\varepsilon\circ\psi),\,\varphi\psi},
	\]
	where $\varphi\psi=\varphi\circ\psi$ and characters are multiplied
	pointwise.
	
	The maps $\varepsilon\mapsto\tau_{\varepsilon,\mathrm{id}}$ and
	$\varphi\mapsto\tau_{1,\varphi}$ are injective group homomorphisms
	into $\Aut_U(U^0)$. Their images have trivial intersection: if
	$\tau_{\varepsilon,\mathrm{id}}=\tau_{1,\varphi}$, then
	$\varepsilon(\lambda)K_\lambda=K_{\varphi(\lambda)}$ for every
	$\lambda\in Q$, so the linear independence of the elements $K_\lambda$
	implies $\varphi=\mathrm{id}$ and $\varepsilon=1$.
	Moreover,
	\[
	\tau_{\varepsilon,\varphi}
	=\tau_{\varepsilon\circ\varphi^{-1},\mathrm{id}}
	\circ\tau_{1,\varphi}.
	\]
	Since $\varepsilon\mapsto\varepsilon\circ\varphi^{-1}$ is a bijection
	of $\mathcal E$, the set of all $\tau_{\varepsilon,\varphi}$ is the
	product of these two subgroups.
	
	For $\lambda\in Q$, direct evaluation gives
	\[
	\bigl(\tau_{1,\varphi}\circ\tau_{\varepsilon,\mathrm{id}}
	\circ\tau_{1,\varphi}^{-1}\bigr)(K_\lambda)
	=\varepsilon\bigl(\varphi^{-1}(\lambda)\bigr)K_\lambda
	=\tau_{\varepsilon\circ\varphi^{-1},\mathrm{id}}(K_\lambda).
	\]
	Thus the second subgroup normalizes the first, with the stated action.
	Their product is therefore a subgroup of $\Aut_U(U^0)$ and is the
	internal semidirect product of these two subgroups.
	Explicitly, the resulting isomorphism is given by
	\[
	\tau_{\varepsilon,\varphi}
	\longmapsto (\varepsilon\circ\varphi^{-1},\varphi).
	\qedhere
	\]
\end{proof}

\subsection{Twisted modules and type-1 representations}\label{subsec:twisted-modules}

For a $U$-module $M$ and $\mu\in P$, set
\[
M_\mu
=
\{v\in M:
K_i v=q_i^{\langle\mu,\alpha_i^\vee\rangle}v
\text{ for all }i\}.
\]
We call $M$ a type-1 module if
\[
M=\bigoplus_{\mu\in P}M_\mu.
\]
\begin{definition}\label{def:twisted-module}
	Let $M$ be a left $U$-module and let $\sigma\in\Aut(U)$.
	The $\sigma$-twisted module $M^\sigma$ is the vector space $M$
	equipped with the $U$-action
	\[
	x\cdot_\sigma m=\sigma(x)m
	\qquad(x\in U,\ m\in M).
	\]
\end{definition}

Since $\sigma(U)=U$, a subspace $N\subseteq M$ is a submodule of
$M^\sigma$ if and only if it is a submodule of $M$.
Hence $M^\sigma$ is irreducible if and only if $M$ is irreducible. We next record the effect of twisting on weights.

\begin{lemma}\label{lem:weight-transform}
	Let $\sigma\in\Aut(U)$ satisfy
	\[
	\sigma(K_i)=K_1^{b_{i1}}\cdots K_n^{b_{in}}
	\qquad(1\le i\le n).
	\]
	Let $M$ be a type-1 $U$-module and let $v\in M$ be a nonzero weight
	vector of weight $\mu\in P$. If $v$, viewed as a vector of $M^\sigma$,
	has a type-1 weight $\mu'\in P$, then
	\[
	[\mu']=D^{-1}A_\sigma D[\mu].
	\]
\end{lemma}

\begin{proof}
	Write
	\[
	\mu=\sum_{j=1}^{n}r_j\omega_j,
	\qquad
	\mu'=\sum_{j=1}^{n}r_j'\omega_j,
	\qquad
	r_j,r_j'\in\mathbb Z.
	\]
	Since $K_jv=q_j^{r_j}v=q^{d_jr_j}v$, Definition~2.8 gives
	\[
		K_i\cdot_\sigma v
		=\sigma(K_i)v
		=\left(\prod_{j=1}^{n}q_j^{b_{ij}r_j}\right)v
		=q^{\sum_{j=1}^{n}b_{ij}d_jr_j}v.
	\]
	On the other hand, the assumption on $\mu'$ gives
	\[
	K_i\cdot_\sigma v=q_i^{r_i'}v=q^{d_ir_i'}v.
	\]
	Since $v\ne0$ and $q$ is not a root of unity, comparison yields
	\[
	d_ir_i'=\sum_{j=1}^{n}b_{ij}d_jr_j
	\qquad(1\le i\le n).
	\]
	These equalities are equivalent to
	\[
	D[\mu']=A_\sigma D[\mu].
	\]
	Multiplying by $D^{-1}$ gives the required formula.
\end{proof}
We recall the standard classification of finite-dimensional
irreducible $U$-modules; see~[9, \S\S5.1--5.2 and Theorem~5.10].

\begin{lemma}\label{lem:sign-twist}
Let $V$ be a finite-dimensional irreducible $U$-module.
Then there exist $\lambda\in P^+$ and $\varepsilon\in\mathcal E$
such that
\[
V\cong L(\lambda)^{s_\varepsilon},
\]
where $L(\lambda)$ denotes the finite-dimensional irreducible
type-1 $U$-module of highest weight $\lambda$.
\end{lemma}

\begin{lemma}\label{lem:signfree-twist-type1}
	Let $\tau\in\operatorname{Aut}(U)$ be sign-free on the Cartan torus.
	Then, for every $\lambda\in P^+$, there exists a unique
	$\mu\in P^+$ such that
	\[
	L(\lambda)^\tau\cong L(\mu).
	\]
\end{lemma}

\begin{proof}
	Twisting preserves finite dimensionality and irreducibility.
	By Lemma~2.10, there exist $\mu\in P^+$ and
	$\varepsilon\in\mathcal E$ such that
	\[
	L(\lambda)^\tau\cong L(\mu)^{s_\varepsilon}.
	\]
Since $L(\lambda)$ is type-1, the $K_j$ are simultaneously
	diagonalizable on $L(\lambda)$, with eigenvalues in
	$q_j^{\mathbb Z}\subseteq q^{\mathbb Z}$.
	As $\tau$ is sign-free, each $\tau(K_i)$ is a Laurent monomial
	in the $K_j$. Hence every $K_i$-eigenvalue on $L(\lambda)^\tau$
	belongs to $q^{\mathbb Z}$.
	
	Let $v_\mu$ be a highest-weight vector of $L(\mu)$.
	In $L(\mu)^{s_\varepsilon}$, we have
	\[
	K_i\cdot_{s_\varepsilon}v_\mu
	=\varepsilon(\alpha_i)
	q_i^{\langle\mu,\alpha_i^\vee\rangle}v_\mu.
	\]
	Thus
	\[
	\varepsilon(\alpha_i)
	q_i^{\langle\mu,\alpha_i^\vee\rangle}\in q^{\mathbb Z}.
	\]
	Since $q_i^{\langle\mu,\alpha_i^\vee\rangle}\in q^{\mathbb Z}$,
	it follows that
	\[
	\varepsilon(\alpha_i)\in
	q^{\mathbb Z}\cap\{\pm1\}=\{1\},
	\]
	where the last equality follows from the assumptions that
	$q$ is not a root of unity and $\operatorname{char}\Bbbk=0$.
	The simple roots generate $Q$, so $\varepsilon=1$.
	Consequently,
	\[
	L(\lambda)^\tau\cong L(\mu).
	\]
	The uniqueness of $\mu$ follows from the fact that irreducible
	highest-weight modules with distinct highest weights are nonisomorphic.
\end{proof}

\section{Preservation of the root system}\label{sec:uniform-matrix}
In this section, we prove that, for a sign-free automorphism $\tau$,
the weight transformation defined by $D^{-1}A_\tau D$ preserves
the root system. We first show that this transformation preserves
the weight lattice and the Weyl-invariant character ring.
We then use the associated fixed field to prove that it normalizes
the Weyl group, and recover the roots from the reflection action
on the weight lattice. This yields the matrix decomposition
in Theorem~3.6.

\subsection{Formal characters and preservation of the character ring}
Denote by
\[
\mathbb Z[P]=\bigoplus_{\mu\in P}\mathbb Z e^\mu
\]
the integral group algebra of $P$, with multiplication
$e^\mu e^\nu=e^{\mu+\nu}$.
The Weyl group acts on $\mathbb Z[P]$ through its action on the weight lattice:
\[
w(e^\mu)=e^{w\mu}\qquad(w\in W,\ \mu\in P).
\]
We denote the corresponding invariant subring by
\[
\mathbb Z[P]^W
:=\{f\in\mathbb Z[P]\mid w(f)=f\text{ for all }w\in W\}.
\]

If $M$ is a finite-dimensional type-$1$ $U$-module, its formal character is
\[
\operatorname{ch}M=\sum_{\mu\in P}(\dim M_\mu)e^\mu\in\mathbb Z[P].
\]
For $\lambda\in P^+$, the weight multiplicities of $L(\lambda)$ are $W$-invariant, and hence
\[
\operatorname{ch}L(\lambda)\in\mathbb Z[P]^W.
\]
For $\lambda\in P^+$, define the Weyl-orbit sum
\[
m_\lambda:=\sum_{\mu\in W\lambda}e^\mu.
\]
Every Weyl orbit in $P$ contains a unique dominant weight,
and the coefficients of a $W$-invariant element are constant
on each orbit. Hence
\[
\mathbb Z[P]^W
=\bigoplus_{\lambda\in P^+}\mathbb Z m_\lambda.
\]

The irreducible characters form a $\mathbb Z$-basis of
$\mathbb Z[P]^W$. Indeed, each weight in the Weyl orbit $W\lambda$
occurs in $L(\lambda)$ with multiplicity one. Every other Weyl
orbit occurring in the character has a unique dominant
representative $\mu$ satisfying $\mu<\lambda$ in the dominance
order. Since the character is $W$-invariant, grouping its terms
by Weyl orbits gives
\[
\operatorname{ch}L(\lambda)=m_\lambda+
\sum_{\substack{\mu\in P^+\\\mu<\lambda}}
a_{\lambda\mu}m_\mu,
\qquad a_{\lambda\mu}\in\mathbb Z_{\geq0}.
\]
For each $\lambda\in P^+$, there are only finitely many
$\mu\in P^+$ with $\mu\leq\lambda$. These unitriangular
relations therefore give
\[
\mathbb Z[P]^W
=\bigoplus_{\lambda\in P^+}\mathbb Z\,\operatorname{ch}L(\lambda).
\]
We call $\mathbb Z[P]^W$ the Weyl-invariant character ring.

Let $\tau\in\Aut(U)$ be sign-free on $U^0$ and put
\[
Y:=D^{-1}A_\tau D\in \mathrm{GL}_n(\mathbb Q).
\]
Since $\omega_1,\ldots,\omega_n$ form a basis of $P\otimes_{\mathbb Z}\mathbb Q$, the matrix $Y$ determines a unique $\mathbb Q$-linear automorphism
\[
\psi_Y:P\otimes_{\mathbb Z}\mathbb Q\longrightarrow P\otimes_{\mathbb Z}\mathbb Q
\]
characterized by
\begin{equation}\label{eq:psiY-definition}
[\psi_Y(\mu)]=Y[\mu]
\qquad(\mu\in P\otimes_{\mathbb Z}\mathbb Q).
\end{equation}

\begin{lemma}\label{lem:Y-preserves-P}
	The linear map $\psi_Y$ preserves the weight lattice:
	\[
	\psi_Y(P)=P.
	\]
	In particular, $Y\in \mathrm{GL}_n(\mathbb Z)$.
\end{lemma}

\begin{proof}
	For each fundamental weight $\omega_i$, choose a nonzero highest-weight
	vector $v_i\in L(\omega_i)$.
	Since each $\tau(K_j)$ is a Laurent monomial in $K_1,\ldots,K_n$, the
	vector $v_i$ is also a common eigenvector for the $K_j$ in the twisted
	module $L(\omega_i)^\tau$.
	By Lemma~2.11, this module is type-1, so the weight of $v_i$ in the
	twisted module belongs to $P$.
    By Lemma~2.9, the coordinate column of this weight is
$Y[\omega_i]=[\psi_Y(\omega_i)]$.
Hence the weight is $\psi_Y(\omega_i)$. Thus
	\[
	\psi_Y(\omega_i)\in P
	\qquad (1\leq i\leq n).
	\]
Since $\omega_1,\ldots,\omega_n$ form a $\mathbb Z$-basis
of $P$, the matrix $Y$ has integer entries. Moreover,
\[
\det Y
=\det(D^{-1}A_\tau D)
=\det A_\tau
=\pm1.
\]
The adjugate formula shows that $Y^{-1}$ also has integer
entries. Hence $Y\in\mathrm{GL}_n(\mathbb Z)$ and
\mbox{$\psi_Y(P)=P$.}
\end{proof}

By Lemma~\ref{lem:Y-preserves-P}, $\psi_Y$ induces a $\mathbb Z$-algebra automorphism
\[
\Psi_Y:\mathbb Z[P]\longrightarrow\mathbb Z[P],
\qquad
\Psi_Y(e^\mu)=e^{\psi_Y(\mu)}.
\]
\begin{lemma}\label{lem:character-ring-preserved}
	The automorphism $\Psi_Y$ preserves the Weyl-invariant character ring:
	\[
	\Psi_Y\bigl(\mathbb Z[P]^W\bigr)=\mathbb Z[P]^W.
	\]
	More precisely, $\Psi_Y$ permutes the irreducible characters
	$\operatorname{ch}L(\lambda)$, $\lambda\in P^+$.
\end{lemma}

\begin{proof}
	Let $\lambda\in P^+$. By Lemma~2.11, there is a unique
	$\lambda'\in P^+$ such that
	\[
	L(\lambda)^\tau\cong L(\lambda').
	\]
	By Lemma~2.9, twisting sends a weight $\mu$ to $\psi_Y(\mu)$.
	Since $\psi_Y$ is bijective on $P$ by Lemma~3.1, we have
	\[
	\bigl(L(\lambda)^\tau\bigr)_{\psi_Y(\mu)}
	=L(\lambda)_\mu
	\qquad(\mu\in P)
	\]
	as subspaces of the common underlying vector space. Therefore
	\[
		\Psi_Y\bigl(\operatorname{ch}L(\lambda)\bigr)
		=\sum_{\mu\in P}
		\dim_{\Bbbk}L(\lambda)_\mu\,e^{\psi_Y(\mu)}
		=\operatorname{ch}\bigl(L(\lambda)^\tau\bigr)
		=\operatorname{ch}L(\lambda').
	\]
	The automorphism $\tau^{-1}$ is also sign-free. Applying
	Lemma~2.11 to $\tau^{-1}$ and using that twisting by $\tau$
	and by $\tau^{-1}$ are mutually inverse shows that
	$\lambda\mapsto\lambda'$ is bijective on $P^+$.
	Thus $\Psi_Y$ permutes the irreducible characters.
	
	Since the irreducible characters form a $\mathbb Z$-basis of
	$\mathbb Z[P]^W$, it follows that
	\[
	\Psi_Y\bigl(\mathbb Z[P]^W\bigr)=
	\bigoplus_{\lambda\in P^+}
	\mathbb Z\,\Psi_Y\bigl(\operatorname{ch}L(\lambda)\bigr)=
	\bigoplus_{\lambda'\in P^+}
	\mathbb Z\,\operatorname{ch}L(\lambda')
	=\mathbb Z[P]^W.
	\qedhere
	\]
\end{proof}

\subsection{The fixed field and normalization of the Weyl group}

Let
\[
F:=\operatorname{Frac}\bigl(\mathbb Q[P]\bigr).
\]
The actions of $W$ and $\Psi_Y$ on $\mathbb Z[P]$ extend naturally
to $\mathbb Q[P]$ and $F$. We write $\mathbb Q[P]^W$ for the
invariant subring of $\mathbb Q[P]$ and $F^W$ for the fixed field
of $W$ in $F$.

\begin{lemma}\label{lem:fixed-field-frac}
	One has
	\[
	F^W=\operatorname{Frac}\bigl(\mathbb Q[P]^W\bigr).
	\]
\end{lemma}

\begin{proof}
	The inclusion $\operatorname{Frac}(\mathbb Q[P]^W)\subseteq F^W$
	is immediate. Conversely, let $f=a/b\in F^W$, where
	$a,b\in\mathbb Q[P]$ and $b\ne0$. Set
	\[
	B:=\prod_{w\in W}w(b)\in\mathbb Q[P]^W\setminus\{0\}.
	\]
	Since $b$ divides $B$, we have $fB=a(B/b)\in\mathbb Q[P]$.
	Both $f$ and $B$ are $W$-invariant, so $fB\in\mathbb Q[P]^W$.
	Therefore
	\[
	f=\frac{fB}{B}\in\operatorname{Frac}\bigl(\mathbb Q[P]^W\bigr),
	\]
	which proves the reverse inclusion.
\end{proof}

%
%
%
%
%

\begin{lemma}\label{lem:Y-normalizes-W}
	The automorphism $\psi_Y$ normalizes $W$ in
	$\operatorname{Aut}_{\mathbb Z}(P)$; that is,
	\[
	\psi_YW\psi_Y^{-1}=W.
	\]
\end{lemma}

\begin{proof}
	By Lemma~3.2,
	\[
	\Psi_Y\bigl(\mathbb Z[P]^W\bigr)=\mathbb Z[P]^W.
	\]
	Every element of $\mathbb Q[P]^W$ has a positive integer multiple
	in $\mathbb Z[P]^W$. Extending $\Psi_Y$ $\mathbb Q$-linearly
	therefore gives
	\[
	\Psi_Y\bigl(\mathbb Q[P]^W\bigr)=\mathbb Q[P]^W.
	\]
	Since $F^W=\operatorname{Frac}(\mathbb Q[P]^W)$ by Lemma~3.3,
	it follows that
	\[
	\Psi_Y(F^W)=F^W.
	\]
	
	The action of $W$ on $F$ is faithful.
	Indeed, if $w\in W$ acts trivially on $F$, then
	\[
	e^{w\mu}=w(e^\mu)=e^\mu
	\qquad(\mu\in P).
	\]
	Since distinct weights label distinct basis elements of
	$\mathbb Q[P]$, we have $w\mu=\mu$ for every $\mu\in P$.
	The faithful action of $W$ on $P$ therefore implies $w=1$.
	Since $W$ is finite, Artin's fixed-field theorem gives
	\[
	\operatorname{Aut}_{F^W}(F)=W.
	\]
	Now let $w\in W$. For every $f\in F^W$, we have
	$\Psi_Y^{-1}(f)\in F^W$, and hence
	\[
	(\Psi_Yw\Psi_Y^{-1})(f)
	=\Psi_Y\bigl(w(\Psi_Y^{-1}(f))\bigr)
	=\Psi_Y(\Psi_Y^{-1}(f))
	=f.
	\]
	Thus $\Psi_Yw\Psi_Y^{-1}$ fixes $F^W$ pointwise, so
	$\Psi_Yw\Psi_Y^{-1}=u$ for some $u\in W$.
	For every $\mu\in P$, this gives
	\[
	e^{(\psi_Yw\psi_Y^{-1})(\mu)}
	=(\Psi_Yw\Psi_Y^{-1})(e^\mu)
	=u(e^\mu)
	=e^{u\mu}.
	\]
	Since the elements $e^\mu$, $\mu\in P$, are distinct basis elements
	of $\mathbb Q[P]$, we obtain $\psi_Yw\psi_Y^{-1}=u\in W$.
	Consequently,
	\[
	\psi_YW\psi_Y^{-1}\subseteq W.
	\]
	Both groups have the same finite order, so this inclusion is an equality.
\end{proof}

\subsection{Root-system automorphisms}
We now use the reflection action of $W$ on $P$ to show that
the lattice automorphism obtained above preserves the root system.

\begin{proposition}\label{prop:root-system-recovery}
	Let $\tau\in\operatorname{Aut}(U)$ be sign-free on the Cartan torus,
	let $Y=D^{-1}A_\tau D$, and let $\psi_Y$ be defined by \eqref{eq:psiY-definition}.
	Then
	\[
	\psi_Y(\Phi)=\Phi.
	\]
\end{proposition}

\begin{proof}
	By Lemmas~3.1 and~3.4,
	\[
	\psi_Y(P)=P,
	\qquad
	\psi_YW\psi_Y^{-1}=W.
	\]
	We first prove that
	\[
	(1-s_\gamma)P=\mathbb Z\gamma
	\qquad(\gamma\in\Phi).
	\]
	For every $\lambda\in P$, the reflection formula gives
	\[
	(1-s_\gamma)(\lambda)
	=\langle\lambda,\gamma^\vee\rangle\gamma
	\in\mathbb Z\gamma.
	\]
	Thus $(1-s_\gamma)P\subseteq\mathbb Z\gamma$.
	
	Conversely, write $\gamma=w(\alpha_i)$ for some $w\in W$
	and some simple root $\alpha_i$. Since $w(\omega_i)\in P$
	and $\gamma^\vee=w(\alpha_i^\vee)$, we have
	\[
	\langle w(\omega_i),\gamma^\vee\rangle
	=\langle w(\omega_i),w(\alpha_i^\vee)\rangle
	=\langle\omega_i,\alpha_i^\vee\rangle
	=1.
	\]
	Consequently, for every $m\in\mathbb Z$,
	\[
	(1-s_\gamma)\bigl(mw(\omega_i)\bigr)=m\gamma.
	\]
	Hence $\mathbb Z\gamma\subseteq(1-s_\gamma)P$,
	which proves the identity.
	
	Now let $\alpha\in\Phi$.
	The conjugate $\psi_Ys_\alpha\psi_Y^{-1}$ belongs to $W$
	and has a fixed subspace of codimension one.
	Since $W$ acts orthogonally, this conjugate is a reflection. Since every reflection in $W$
	is a root reflection, there exists $\beta\in\Phi$ such that
	\[
	\psi_Ys_\alpha\psi_Y^{-1}=s_\beta.
	\]
	Using the identity above and $\psi_Y(P)=P$, we obtain
	\[
		\mathbb Z\psi_Y(\alpha)
		=\psi_Y\bigl((1-s_\alpha)P\bigr)
		=(1-\psi_Ys_\alpha\psi_Y^{-1})P
		=(1-s_\beta)P
		=\mathbb Z\beta.
	\]
	Therefore $\psi_Y(\alpha)=\pm\beta\in\Phi$,
	and hence $\psi_Y(\Phi)\subseteq\Phi$.
	Since $\Phi$ is finite and $\psi_Y$ is injective,
	this inclusion is an equality.
\end{proof}

\begin{theorem}\label{thm:uniform-root-rigidity}
	Let $\tau\in\Aut(U)$ be sign-free on the Cartan torus.
	Then there exist unique elements $\xi\in\Gamma$ and $w\in W$ such that
	\[
	D^{-1}A_\tau D=[\xi]\vartheta(w).
	\]
\end{theorem}

\begin{proof}
	Let $Y=D^{-1}A_\tau D$. By Proposition~3.5,
	$\psi_Y\in\Aut(\Phi)$. Lemma~2.6 gives a unique factorization
	\[
	\psi_Y=w_0\xi,\qquad w_0\in W,\quad \xi\in\Gamma.
	\]
	Since $W$ is normal in $\Aut(\Phi)$, we may set
	$w=\xi^{-1}w_0\xi\in W$. Then
	\[
	\psi_Y=w_0\xi=\xi(\xi^{-1}w_0\xi)=\xi w.
	\]
	The uniqueness of $\xi$ and $w$ follows from Lemma~2.6.
	Taking matrices with respect to the ordered basis
	$\omega_1,\ldots,\omega_n$ gives
	\[
	D^{-1}A_\tau D=Y=[\psi_Y]=[\xi]\vartheta(w).\qedhere
	\]
\end{proof}

\begin{corollary}\label{cor:uniform-matrix}
	For every $\sigma\in\operatorname{Aut}(U)$, there exist unique elements
	$\xi\in\Gamma$ and $w\in W$ such that
	\[
	D^{-1}A_\sigma D=[\xi]\vartheta(w).
	\]
\end{corollary}

\begin{proof}
	By Lemma~2.4, choose $\delta\in\mathcal E$ such that
	$\tau=s_\delta\circ\sigma$ is sign-free on the Cartan torus
	and $A_\tau=A_\sigma$.
	Existence and uniqueness now follow from Theorem~3.6.
\end{proof}

\section{Induced automorphism groups of the Cartan torus and center}
In this section, we determine $\Aut_U(U^0)$ and
$\Aut_U(Z(U))$ using the results of Section~3.
We also prove $\ker\rho_0\subseteq\ker\rho_Z$ and determine
the quotient $\ker\rho_Z/\ker\rho_0$.
\subsection{The induced automorphism group of the Cartan torus}
We now determine $\Aut_U(U^0)$. The main remaining step is to relate
the weight transformation obtained in Section~3 to the automorphism of
the root lattice induced by the toral restriction.

\begin{theorem}\label{prop:cartan-criterion-general}
	A $\Bbbk$-algebra automorphism $\tau$ of $U^0$ is the restriction of
	an algebra automorphism of $U$ if and only if there exist unique elements
	$\varepsilon\in\mathcal{E}$ and $\varphi\in\Aut(\Phi)$ such that
	\[
	\tau(K_\lambda)=\varepsilon(\lambda)K_{\varphi(\lambda)}
	\qquad (\lambda\in Q).
	\]
	Consequently,
	\[
	\Aut_U(U^0)\cong\mathcal{E}\rtimes\Aut(\Phi),
	\]
	where $\varphi\cdot\varepsilon=\varepsilon\circ\varphi^{-1}$.
\end{theorem}

\begin{proof}
	Suppose that $\tau\in\Aut_U(U^0)$, and choose $\sigma\in\Aut(U)$
	with $\sigma|_{U^0}=\tau$. By Lemma~2.3, there exist unique elements
	$\varepsilon\in\mathcal{E}$ and $\varphi\in\Aut_{\mathbb Z}(Q)$ such that
	\[
	\tau(K_\lambda)=\sigma(K_\lambda)
	=\varepsilon(\lambda)K_{\varphi(\lambda)}
	\qquad (\lambda\in Q).
	\]
	Write $A_\sigma=(b_{ij})$ and put $Y=D^{-1}A_\sigma D=(y_{ij})$.
	By the definition of $A_\sigma$,
	\[
	\varphi(\alpha_i)=\sum_{k=1}^{n}b_{ik}\alpha_k
	\qquad (1\leq i\leq n).
	\]
	By Corollary~3.7, there exist
	$\xi\in\Gamma$ and $w\in W$ such that
	\[
	[\psi_Y]=Y=[\xi]\vartheta(w).
	\]
	Hence $\psi_Y=\xi w\in\Aut(\Phi)$.
	
	Since $A_\sigma D=DY$, we have
	\begin{equation}\label{eq:toral-weight-entries}
		b_{ij}d_j=d_i y_{ij}\qquad (1\leq i,j\leq n).
	\end{equation}
	Using $(\alpha_k,\omega_j)=d_k\delta_{kj}$ and
	$\psi_Y(\omega_j)=\sum_{k=1}^{n}y_{kj}\omega_k$, we obtain
	\[
	\begin{aligned}
		(\varphi(\alpha_i),\omega_j)
		&=\sum_{k=1}^{n}b_{ik}(\alpha_k,\omega_j)=b_{ij}d_j,\\
		(\alpha_i,\psi_Y(\omega_j))
		&=\sum_{k=1}^{n}y_{kj}(\alpha_i,\omega_k)=d_i y_{ij}.
	\end{aligned}
	\]
	Since both $W$ and $\Gamma$ preserve the standard bilinear form,
	so does $\psi_Y=\xi w$.
	Together with \eqref{eq:toral-weight-entries}, this gives
	\[
		(\varphi(\alpha_i),\omega_j)
		=(\alpha_i,\psi_Y(\omega_j))
		=(\psi_Y^{-1}(\alpha_i),\omega_j)
		\qquad (1\leq i,j\leq n).
	\]
	The fundamental weights form a basis of the rational weight space, and
	the bilinear form is nondegenerate. Therefore
	\[
	\varphi(\alpha_i)=\psi_Y^{-1}(\alpha_i)
	\qquad (1\leq i\leq n).
	\]
	Thus $\varphi$ agrees with $\psi_Y^{-1}$ on $Q$. Since $\psi_Y^{-1}$
	preserves $\Phi$, we obtain $\varphi\in\Aut(\Phi)$.
	
	Conversely, Proposition~2.7 shows that every map of the stated form
	extends to an algebra automorphism of $U$, and that these maps form a
	group isomorphic to $\mathcal{E}\rtimes\Aut(\Phi)$ with the stated action.
	This proves both the converse and the semidirect-product description.
\end{proof}

\noindent\textbf{Remark.}
Since $U^0=\Bbbk[Q]$ is a Laurent polynomial algebra,
its units are precisely the elements $cK_\lambda$ with
$c\in\Bbbk^\times$ and $\lambda\in Q$.
It follows that
\[
\Aut(U^0)
\cong \Hom(Q,\Bbbk^\times)\rtimes\Aut_{\mathbb Z}(Q)
\cong (\Bbbk^\times)^n\rtimes\mathrm{GL}_n(\mathbb Z),
\]
where $\varphi$ acts on a character $\chi$ by
$\varphi\cdot\chi=\chi\circ\varphi^{-1}$.
Theorem~4.1 identifies the subgroup of automorphisms
extending to $U$ as
\[
\Aut_U(U^0)
\cong \Hom(Q,\{\pm1\})\rtimes\Aut(\Phi).
\]
Thus extendibility to $U$ requires precisely that the
character take values in $\{\pm1\}$ and that the lattice
automorphism preserve $\Phi$.

\begin{corollary}
	For each Cartan type, the group of automorphisms induced on the Cartan
	torus is given by
	\[
	\Aut_U(U^0)\cong
	\begin{cases}
		(\mathbb Z/2\mathbb Z)^2,
		& A_1,\\
		(\mathbb Z/2\mathbb Z)^n\rtimes
		\bigl(S_{n+1}\times\mathbb Z/2\mathbb Z\bigr),
		& A_n,\ n\geq2,\\
		(\mathbb Z/2\mathbb Z)^n\rtimes
		\bigl((\mathbb Z/2\mathbb Z)^n\rtimes S_n\bigr),
		& B_n,\ C_n,\ n\geq2,\\
		(\mathbb Z/2\mathbb Z)^4\rtimes W(F_4),
		& D_4,\\
		(\mathbb Z/2\mathbb Z)^n\rtimes
		\bigl((\mathbb Z/2\mathbb Z)^n\rtimes S_n\bigr),
		& D_n,\ n\geq5,\\
		(\mathbb Z/2\mathbb Z)^6\rtimes
		\bigl(W(E_6)\rtimes\mathbb Z/2\mathbb Z\bigr),
		& E_6,\\
		(\mathbb Z/2\mathbb Z)^7\rtimes W(E_7),
		& E_7,\\
		(\mathbb Z/2\mathbb Z)^8\rtimes W(E_8),
		& E_8,\\
		(\mathbb Z/2\mathbb Z)^4\rtimes W(F_4),
		& F_4,\\
		(\mathbb Z/2\mathbb Z)^2\rtimes D_6,
		& G_2.
	\end{cases}
	\]
	Here $D_6$ denotes the dihedral group of order $12$.
	In each outer semidirect product, the first factor is identified with
	$\mathcal E=\Hom(Q,\{\pm1\})$, where $Q$ is the root lattice of the
	indicated type. The second factor is identified with $\Aut(\Phi)$ as
	in the proof, and acts by
	\[
	(\varphi\cdot\varepsilon)(\lambda)
	=\varepsilon\bigl(\varphi^{-1}(\lambda)\bigr)
	\qquad
	(\varphi\in\Aut(\Phi),\ \varepsilon\in\mathcal E,\ \lambda\in Q).
	\]
\end{corollary}

\begin{proof}
	By Theorem~4.1 and Lemma~2.6,
	\[
	\Aut_U(U^0)\cong\mathcal E\rtimes\Aut(\Phi),
	\qquad
	\Aut(\Phi)=W\rtimes\Gamma.
	\]
	It remains to identify the root-system automorphism group in each type.
	Here $\Aut(A_n)$ denotes the automorphism group of a root system of
	type $A_n$, and similarly for the other Cartan types.
	
	For $A_1$, we have $W=\{1,-\mathrm{id}\}$ and
	$\mathcal E\cong\mathbb Z/2\mathbb Z$. Since
	\[
	\bigl((-\mathrm{id})\cdot\varepsilon\bigr)(\lambda)
	=\varepsilon(-\lambda)
	=\varepsilon(\lambda)^{-1}
	=\varepsilon(\lambda)
	\qquad
	(\varepsilon\in\mathcal E,\ \lambda\in Q),
	\]
	the action of $W$ on $\mathcal E$ is trivial. Thus
	\[
	\Aut_U(U^0)\cong\mathcal E\times W
	\cong(\mathbb Z/2\mathbb Z)^2.
	\]
	
	For $A_n$ with $n\geq2$, we have $W(A_n)\cong S_{n+1}$.
	Moreover, $-\mathrm{id}\in\Aut(A_n)\setminus W(A_n)$
	centralizes $W(A_n)$. Since
	$[\Aut(A_n):W(A_n)]=|\Gamma|=2$, we obtain
	\[
	\Aut(A_n)
	=W(A_n)\times\langle-\mathrm{id}\rangle
	\cong S_{n+1}\times\mathbb Z/2\mathbb Z.
	\]
	
	For $B_n$ and $C_n$, the Dynkin-diagram automorphism group is trivial,
	and the Weyl group is the full group of signed coordinate permutations.
	Hence
	\[
	\Aut(B_n)\cong\Aut(C_n)
	\cong(\mathbb Z/2\mathbb Z)^n\rtimes S_n.
	\]
	
	For $D_n$ with $n\geq5$, use the standard coordinate basis
	$e_1,\ldots,e_n$. Then
	\[
	W(D_n)=
	\left\{
	e_i\mapsto\varepsilon_i e_{\pi(i)}
	\;\middle|\;
	\pi\in S_n,\ \varepsilon_i\in\{\pm1\},\
	\prod_{i=1}^n\varepsilon_i=1
	\right\},
	\]
	and
	\[
	\Gamma=\langle\xi\rangle\cong\mathbb Z/2\mathbb Z,
	\qquad
	\xi(e_i)=
	\begin{cases}
		e_i, & 1\leq i<n,\\
		-e_n, & i=n.
	\end{cases}
	\]
	Hence
	\[
		\Aut(D_n)
		=W(D_n)\rtimes\Gamma
		=
		\left\{
		e_i\mapsto\varepsilon_i e_{\pi(i)}
		\;\middle|\;
		\pi\in S_n,\ \varepsilon_i\in\{\pm1\}
		\right\}
		\cong(\mathbb Z/2\mathbb Z)^n\rtimes S_n.
\]
	For $D_4$, identify its standard root system with the long-root
	subsystem of $F_4$:
	\[
	\Phi(D_4)=\Phi_{\mathrm{long}}(F_4)
	=\{\pm e_i\pm e_j:1\leq i<j\leq4\}.
	\]
	The Weyl group $W(F_4)$ preserves $\Phi(F_4)$ and root lengths,
	so it preserves $\Phi(D_4)$. Since $\Phi(D_4)$ spans the ambient
	space, this action is faithful. Hence
	\[
	W(F_4)\subseteq\Aut(D_4).
	\]
	Since $\Gamma\cong S_3$ in type $D_4$, we have
	\[
	|\Aut(D_4)|
	=|W(D_4)|\,|\Gamma|
	=(2^3\cdot4!)\cdot6
	=1152
	=|W(F_4)|.
	\]
	The inclusion is therefore an equality, giving
	$\Aut(D_4)\cong W(F_4)$.
	
	Finally, $\Gamma\cong\mathbb Z/2\mathbb Z$ in type $E_6$ and is
	trivial in types $E_7$, $E_8$, $F_4$, and $G_2$.
	Also, $W(G_2)\cong D_6$.
	Substituting these descriptions into Theorem~4.1 gives all the stated
	groups, with the specified actions.
\end{proof}

\subsection{Central characters and the sign-free action on the center}
We next use the unique matrix decomposition from Theorem~\ref{thm:uniform-root-rigidity}
to determine how a sign-free automorphism acts on finite-dimensional irreducible
type-$1$ modules. This will be the main input for determining its induced action
on the center.

\begin{lemma}
	Let $\tau\in\Aut(U)$ be sign-free on the Cartan torus, and let
	$\xi\in\Gamma$ and $w\in W$ be the unique elements such that
	\[
	D^{-1}A_\tau D=[\xi]\vartheta(w).
	\]
	Then, for every dominant integral weight $\lambda\in P^+$,
	\[
	L(\lambda)^\tau\cong L(\xi(\lambda)).
	\]
\end{lemma}

\begin{proof}\label{lem:all-twists}
	By Lemma~2.11, $L(\lambda)^\tau$ is a finite-dimensional irreducible
	type-1 module. By Lemma~2.9, a weight $\mu$ of $L(\lambda)$ is
	transformed under twisting into a weight $\mu'$ satisfying
	\[
	[\mu']
	=D^{-1}A_\tau D[\mu]
	=[\xi]\vartheta(w)[\mu]
	=[\xi(w(\mu))].
	\]
	Thus $\mu'=\xi(w(\mu))$. Since $\xi w$ is invertible and the weight
	set of $L(\lambda)$ is $W$-stable, we obtain
	\[
	\operatorname{wt}\bigl(L(\lambda)^\tau\bigr)
	=\xi w\bigl(\operatorname{wt} L(\lambda)\bigr)
	=\xi\bigl(\operatorname{wt} L(\lambda)\bigr).
	\]
	
	Since $\xi$ permutes the simple roots, it preserves $Q^{+}$ and $P^+$.
	For every $\mu\in\operatorname{wt} L(\lambda)$, we have $\lambda-\mu\in Q^{+}$,
	and hence
	\[
	\xi(\lambda)-\xi(\mu)
	=\xi(\lambda-\mu)\in Q^{+}.
	\]
	Moreover, $\lambda\in\operatorname{wt} L(\lambda)$, so $\xi(\lambda)$ occurs as
	a weight of $L(\lambda)^\tau$. Therefore $L(\lambda)^\tau$ has highest
	weight $\xi(\lambda)$.
	The classification by dominant highest weights now gives
	\[
	L(\lambda)^\tau\cong L(\xi(\lambda)).
\qedhere	\]
\end{proof}

For $\lambda\in P^+$, let
\[
\chi_\lambda:Z(U)\longrightarrow\Bbbk
\]
be the central character of $L(\lambda)$. This character is $\Bbbk$-valued
even though $\Bbbk$ is not assumed algebraically closed. Indeed, let
$v_\lambda$ span the one-dimensional highest-weight space of $L(\lambda)$.
Every $z\in Z(U)$ preserves that space, so $zv_\lambda=cv_\lambda$ for
some $c\in\Bbbk$. By centrality,
\[
z(uv_\lambda)=u(zv_\lambda)=c\,uv_\lambda
\qquad(u\in U).
\]
Since $L(\lambda)=Uv_\lambda$, the element $z$ acts on all of
$L(\lambda)$ by the scalar $c=\chi_\lambda(z)$.
The following lemma shows that these central characters separate the
elements of $Z(U)$.

\begin{lemma}\label{lem:separate}
	Let $z\in Z(U)$. If $\chi_\lambda(z)=0$ for all $\lambda\in P^+$,
	then $z=0$.
\end{lemma}

\begin{proof}
	Suppose that $z\ne0$. By Proposition~2.1, the restriction of $\pi$ to
	$Z(U)$ is injective, so $\pi(z)\ne0$. Write
	\[
	\pi(z)=\sum_{j=1}^{r}a_jK_{\nu_j},
	\qquad a_j\in\Bbbk\setminus\{0\},\quad \nu_j\in Q,
	\]
	where $r\geq1$ and the $\nu_j$ are pairwise distinct.
	
	For $\lambda\in P^+$, let $v_\lambda$ be as above.
	By the definition of $\pi$, the element
	$z-\pi(z)$ is a sum of terms in $U^-_{-\nu}U^0U^+_\nu$ with
	$\nu\in Q^{+}\setminus\{0\}$, and hence annihilates $v_\lambda$.
	Since $K_\nu v_\lambda=q^{(\lambda,\nu)}v_\lambda$, we obtain
	\[
	0=\chi_\lambda(z)
	=\sum_{j=1}^{r}a_jq^{(\lambda,\nu_j)}
	\qquad(\lambda\in P^+).
	\]
	
	We choose a strictly dominant integral weight $\eta$ such that
	$(\eta,\nu_1),\ldots,(\eta,\nu_r)$ are pairwise distinct integers.
	For this purpose, consider
	\[
	\eta=\omega_1+t\omega_2+\cdots+t^{n-1}\omega_n,
	\qquad t\in\mathbb Z_{>0}.
	\]
	For each $j\ne k$, the expression $(\eta,\nu_j-\nu_k)$ is a nonzero
	polynomial in $t$, since the fundamental weights form a basis and the
	bilinear form is nondegenerate. There are only finitely many such
	polynomials, each with finitely many roots. Thus a positive integer $t$
	can be chosen so that they are all nonzero, giving the required $\eta$.
	
	Taking $\lambda=m\eta$ for $m=0,\ldots,r-1$ yields
	\[
	\sum_{j=1}^{r}a_j\bigl(q^{(\eta,\nu_j)}\bigr)^m=0
	\qquad(0\leq m\leq r-1).
	\]
	The coefficient matrix is a Vandermonde matrix with determinant
	\[
	\prod_{1\leq j<k\leq r}
	\bigl(q^{(\eta,\nu_k)}-q^{(\eta,\nu_j)}\bigr)\ne0,
	\]
	because the integer exponents are distinct and $q$ is not a root of
	unity. Hence $a_1=\cdots=a_r=0$, contradicting the choice of the
	$a_j$. Therefore $z=0$.
\end{proof}

\begin{proposition}\label{prop:signfree-center-general}
	Let $\tau\in\Aut(U)$ be sign-free on the Cartan torus, and let
	$\xi\in\Gamma$ be the unique Dynkin-diagram automorphism associated
	with $\tau$ by Theorem~3.6. Then
	\[
	\tau|_{Z(U)}=G_{\xi^{-1}}|_{Z(U)}.
	\]
\end{proposition}

\begin{proof}
	By Lemma~4.3,
	\[
	L(\lambda)^\tau\cong L(\xi(\lambda))
	\qquad(\lambda\in P^+).
	\]
	For $z\in Z(U)$, its action on $L(\lambda)^\tau$ is the action of
	$\tau(z)$ on $L(\lambda)$. Therefore,
	\[
	\chi_\lambda(\tau(z))=\chi_{\xi(\lambda)}(z).
	\]
	
	On the other hand, let $v_\lambda$ be as above.
	By the definition of $G_{\xi^{-1}}$,
	\[
	E_i\cdot_{G_{\xi^{-1}}}v_\lambda
	=E_{\xi^{-1}(i)}v_\lambda=0,
	\]
	and
	\[
		K_i\cdot_{G_{\xi^{-1}}}v_\lambda
		=K_{\xi^{-1}(i)}v_\lambda
		=q_{\xi^{-1}(i)}^{\langle\lambda,\alpha_{\xi^{-1}(i)}^\vee\rangle}
		v_\lambda
		=q_i^{\langle\xi(\lambda),\alpha_i^\vee\rangle}v_\lambda.
	\]
	Here we have used
	\[
	q_{\xi^{-1}(i)}=q_i,
	\qquad
	\langle\lambda,\alpha_{\xi^{-1}(i)}^\vee\rangle
	=\langle\xi(\lambda),\alpha_i^\vee\rangle.
	\]
	Thus $v_\lambda$ is a highest-weight vector of
	$L(\lambda)^{G_{\xi^{-1}}}$ of highest weight $\xi(\lambda)$.
	Since twisting preserves irreducibility, it follows that
	\[
	L(\lambda)^{G_{\xi^{-1}}}\cong L(\xi(\lambda)).
	\]
	Consequently,
	\[
	\chi_\lambda(G_{\xi^{-1}}(z))
	=\chi_{\xi(\lambda)}(z)
	=\chi_\lambda(\tau(z))
	\qquad(\lambda\in P^+).
	\]
	Both $\tau(z)$ and $G_{\xi^{-1}}(z)$ belong to $Z(U)$,
	so Lemma~4.4 gives
	\[
	\tau(z)=G_{\xi^{-1}}(z).
	\]
	Since $z\in Z(U)$ was arbitrary, the assertion follows.
\end{proof}

We now show that the restriction of an automorphism to the Cartan torus
determines its restriction to the center.

\begin{corollary}\label{4.6}
	One has
	\[
	\ker\rho_0\subseteq\ker\rho_Z.
	\]
	Consequently, the assignment
	\[
	\Aut_U(U^0)\longrightarrow\Aut_U(Z(U)),
	\qquad
	\sigma|_{U^0}\longmapsto\sigma|_{Z(U)}
	\]
	defines a surjective group homomorphism.
\end{corollary}

\begin{proof}
	Let $\sigma\in\ker\rho_0$. Then $\sigma$ is sign-free on the
	Cartan torus and $A_\sigma=I_n$, so
	\[
	D^{-1}A_\sigma D=I_n.
	\]
	By the uniqueness in Theorem~3.6, the Dynkin-diagram
	automorphism associated with $\sigma$ is the identity.
	Proposition~4.5 therefore gives
	\[
	\sigma|_{Z(U)}=\mathrm{id}_{Z(U)}.
	\]
	Hence $\ker\rho_0\subseteq\ker\rho_Z$.
	
	If $\sigma,\tau\in\Aut(U)$ satisfy
	$\sigma|_{U^0}=\tau|_{U^0}$, then
	\[
	\tau^{-1}\circ\sigma\in\ker\rho_0\subseteq\ker\rho_Z,
	\]
	and consequently
	\[
	\sigma|_{Z(U)}=\tau|_{Z(U)}.
	\]
	Thus the assignment is well defined.
	Compatibility with composition makes this assignment a group
	homomorphism. It is surjective by the definition of
	$\Aut_U(Z(U))$.
\end{proof}

\subsection{The lattice $Q\cap 2P$ and sign automorphisms}
To describe the action of sign automorphisms on $Z(U)$,
we first examine the quotient $(Q\cap2P)/2Q$.

\begin{proposition}\label{prop:lattice-quotient}
	The map
	\[
	\eta:(P/Q)[2]\longrightarrow(Q\cap2P)/2Q,
	\qquad \mu+Q\longmapsto2\mu+2Q
	\]
	is a group isomorphism, where
	\[
	(P/Q)[2]=\{x\in P/Q:2x=0\}
	\]
	is the $2$-torsion subgroup of $P/Q$.
	For each Cartan type, the quotient is given by
	\[
	\begin{array}{c|c}
		\text{type} & (Q\cap2P)/2Q \\
		\hline
		\rule{0pt}{16pt}
		A_n &
		\begin{cases}
			0, & n\text{ even},\\
			\mathbb Z/2\mathbb Z, & n\text{ odd}
		\end{cases}
		\\[6pt]
		B_n,\ C_n,\ E_7 & \mathbb Z/2\mathbb Z
		\\[5pt]
		D_n &
		\begin{cases}
			(\mathbb Z/2\mathbb Z)^2, & n\text{ even},\\
			\mathbb Z/2\mathbb Z, & n\text{ odd}
		\end{cases}
		\\[6pt]
		E_6,\ E_8,\ F_4,\ G_2 & 0
	\end{array}
	\]
\end{proposition}

\begin{proof}
	For $\mu+Q\in(P/Q)[2]$, we have $2\mu\in Q\cap2P$.
	If $\mu+Q=\mu'+Q$, then
	\[
	2\mu-2\mu'=2(\mu-\mu')\in2Q,
	\]
	so $2\mu+2Q=2\mu'+2Q$.
	Thus $\eta$ is well defined and is a group homomorphism
	by additivity. Its inverse is given by
	\[
	\eta^{-1}(\lambda+2Q)=\frac{\lambda}{2}+Q
	\qquad(\lambda\in Q\cap2P).
	\]
	Here $\lambda/2\in P$ and
	$2(\lambda/2+Q)=\lambda+Q=Q$.
	For $\nu\in Q$,
	\[
	\frac{\lambda+2\nu}{2}+Q
	=\frac{\lambda}{2}+\nu+Q
	=\frac{\lambda}{2}+Q.
	\]
	Thus the inverse map is well defined.
	The two formulas are mutually inverse, proving the isomorphism.
	
	For an irreducible root system, the standard quotients $P/Q$
	are as follows; see~[5, Plates~I--IX, (VIII)].
	\[
	\begin{array}{c|c}
		\text{type} & P/Q \\
		\hline
		\rule{0pt}{16pt}
		A_n & \mathbb Z/(n+1)\mathbb Z
		\\[5pt]
		B_n,\ C_n,\ E_7 & \mathbb Z/2\mathbb Z
		\\[5pt]
		D_n &
		\begin{cases}
			(\mathbb Z/2\mathbb Z)^2, & n\text{ even},\\
			\mathbb Z/4\mathbb Z, & n\text{ odd}
		\end{cases}
		\\[6pt]
		E_6 & \mathbb Z/3\mathbb Z
		\\[5pt]
		E_8,\ F_4,\ G_2 & 0
	\end{array}
	\]
	Taking the $2$-torsion subgroup gives the table in the statement.
\end{proof}

\begin{lemma}\label{prop:sign-HC-general}
	Let $\delta\in\mathcal E$. Then, for every $z\in Z(U)$,
	\[
	(\gamma_{-\rho}\circ\pi)(s_\delta(z))
	=s_\delta\bigl((\gamma_{-\rho}\circ\pi)(z)\bigr).
	\]
\end{lemma}

\begin{proof}
	The sign automorphism $s_\delta$ preserves the $Q$-grading and the
	triangular decomposition. Hence, for every $\nu\in Q^{+}$,
	\[
	s_\delta(U^+_\nu)=U^+_\nu,\qquad
	s_\delta(U^-_{-\nu})=U^-_{-\nu},\qquad
	s_\delta(U^0)=U^0.
	\]
	Therefore it preserves the two summands in
	\[
	U_0
	=U^0\oplus
	\bigoplus_{\nu\in Q^{+}\setminus\{0\}}
	U^-_{-\nu}U^0U^+_\nu.
	\]
	Since $\pi$ is the projection onto $U^0$ along the second summand,
	we obtain
	\[
	\pi\circ s_\delta=s_\delta\circ\pi
	\]
on $U_0$. Moreover, for every $\lambda\in Q$,
	\[
		\gamma_{-\rho}(s_\delta(K_\lambda))
		=\delta(\lambda)q^{-(\rho,\lambda)}K_\lambda
		=s_\delta(\gamma_{-\rho}(K_\lambda)).
	\]
	Hence $s_\delta$ and $\gamma_{-\rho}$ commute on $U^0$.
	Combining the two commutation relations gives the result.
\end{proof}

\begin{proposition}\label{prop:sign-action-general}
	Let $\delta\in\mathcal E$. Then
	\[
	s_\delta|_{Z(U)}=\operatorname{id}_{Z(U)}
	\quad\text{if and only if}\quad
	\delta|_{Q\cap2P}=1.
	\]
	The assignment
	\[
	s_\delta|_{Z(U)}
	\longmapsto
	\bigl(\mu+Q\longmapsto\delta(2\mu)\bigr)
	\qquad (\mu+Q\in(P/Q)[2])
	\]
	defines a group isomorphism
	\[
	\{s_\delta|_{Z(U)}:\delta\in\mathcal E\}
	\cong
	\Hom\bigl((P/Q)[2],\{\pm1\}\bigr).
	\]
\end{proposition}

\begin{proof}
	By Proposition~2.1 and Lemma~4.8, the Harish--Chandra
	isomorphism
	\[
	\gamma_{-\rho}\circ\pi:
	Z(U)\xrightarrow{\sim}(U^0_{\mathrm{ev}})^W
	\]
	intertwines the sign actions. Thus $s_\delta$ fixes $Z(U)$ pointwise if and only if it fixes
	$(U^0_{\mathrm{ev}})^W$ pointwise.
		Every $\delta\in\mathcal E$ is trivial on $2Q$, since
	$\delta(2\nu)=\delta(\nu)^2=1$ for $\nu\in Q$.
	For $\lambda\in Q\cap2P$, write $\lambda=2\mu$ with $\mu\in P$.
	Since $w\mu-\mu\in Q$ for every $w\in W$, we have
	\[
	w\lambda-\lambda=2(w\mu-\mu)\in2Q,
	\qquad \delta(w\lambda)=\delta(\lambda).
	\]
	Therefore, the Weyl-orbit sum
	\[
	M_\lambda:=\sum_{\nu\in W\lambda}K_\nu\in(U^0_{\mathrm{ev}})^W
	\]
	satisfies
	\[
	s_\delta(M_\lambda)
	=\sum_{\nu\in W\lambda}\delta(\nu)K_\nu
	=\delta(\lambda)M_\lambda.
	\]
	Since the coefficients of a $W$-invariant element are
	constant on each Weyl orbit, the elements
	\[
	\{M_\lambda:\lambda\in Q\cap2P^+\}
	\]
	form a $\Bbbk$-basis of $(U^0_{\mathrm{ev}})^W$.
	Together with the Harish--Chandra isomorphism and the
	$W$-invariance of $\delta|_{Q\cap2P}$, this proves
	\[
	s_\delta|_{Z(U)}=\operatorname{id}_{Z(U)}
	\quad\text{if and only if}\quad
	\delta|_{Q\cap2P}=1.
	\]
	
	Since every sign character is trivial on $2Q$,
	restriction defines a group homomorphism
	\[
	\mathcal E
	\longrightarrow
	\Hom\bigl((Q\cap2P)/2Q,\{\pm1\}\bigr),
	\qquad
	\delta\longmapsto
	\bigl(\lambda+2Q\longmapsto\delta(\lambda)\bigr).
	\]
	To prove surjectivity, let
	$\chi\in\Hom((Q\cap2P)/2Q,\{\pm1\})$.
	Choose a basis $\bar\lambda_1,\ldots,\bar\lambda_r$
	of $(Q\cap2P)/2Q$ over $\mathbb Z/2\mathbb Z$,
	and extend it to a basis
	$\bar\lambda_1,\ldots,\bar\lambda_n$ of $Q/2Q$.
	The formula
	\[
	\widetilde\chi\left(\sum_{i=1}^{n}a_i\bar\lambda_i\right)
	:=
	\prod_{i=1}^{r}\chi(\bar\lambda_i)^{a_i},
	\qquad a_i\in\{0,1\},
	\]
	defines a character of $Q/2Q$ extending $\chi$.
	Thus
	\[
	\delta(\lambda):=\widetilde\chi(\lambda+2Q)
	\qquad(\lambda\in Q)
	\]
	belongs to $\mathcal E$ and maps to $\chi$.
	
	The restriction homomorphism and the homomorphism
	$\delta\mapsto s_\delta|_{Z(U)}$ have the same kernel,
	namely
	\[
	\{\delta\in\mathcal E:\delta|_{Q\cap2P}=1\}.
	\]
	The first isomorphism theorem therefore gives
	\[
	\{s_\delta|_{Z(U)}:\delta\in\mathcal E\}
	\cong
	\Hom\bigl((Q\cap2P)/2Q,\{\pm1\}\bigr),
	\]
	with
	\[
	s_\delta|_{Z(U)}
	\longmapsto
	\bigl(\lambda+2Q\longmapsto\delta(\lambda)\bigr).
	\]
	Finally, composing characters with the isomorphism
	\[
	\eta:(P/Q)[2]\xrightarrow{\sim}(Q\cap2P)/2Q,
	\qquad
	\mu+Q\longmapsto2\mu+2Q,
	\]
	from Proposition~4.7 gives the asserted isomorphism
	\[
	s_\delta|_{Z(U)}
	\longmapsto
	\bigl(\mu+Q\longmapsto\delta(2\mu)\bigr).
	\qedhere
	\]
\end{proof}

\subsection{The induced automorphism group of the center}

We now combine Propositions~4.5 and~4.9 to characterize the automorphisms of the center that are restrictions of algebra automorphisms of $U$ and to determine their group structure.

\begin{theorem}\label{prop:center-criterion-general}
	A $\Bbbk$-algebra automorphism $\tau$ of $Z(U)$ is the restriction of an algebra automorphism of $U$ if and only if there exist $\delta\in\mathcal E$ and $\xi\in\Gamma$ such that
	\[
	\tau=(s_\delta\circ G_{\xi^{-1}})|_{Z(U)}.
	\]
	Moreover, $\xi$ and $\delta|_{Q\cap2P}$ are uniquely determined by $\tau$. Consequently,
	\[
	\Aut_U(Z(U))\cong\Hom\bigl((P/Q)[2],\{\pm1\}\bigr)\rtimes\Gamma,
	\]
	where the action of $\Gamma$ is given by
	\[
	(\xi\cdot\chi)(x)=\chi(\xi^{-1}x)
	\qquad
	\bigl(\xi\in\Gamma,\ \chi\in\Hom((P/Q)[2],\{\pm1\}),\ x\in(P/Q)[2]\bigr).
	\]
\end{theorem}

\begin{proof}
	Suppose that $\tau=\sigma|_{Z(U)}$ for some $\sigma\in\Aut(U)$.
	By Lemma~2.4, choose $\delta\in\mathcal E$ such that
	$s_\delta\circ\sigma$ is sign-free on the Cartan torus, and let
	$\xi\in\Gamma$ be its associated Dynkin-diagram automorphism
	from Theorem~3.6. Proposition~4.5 and $s_\delta^2=\operatorname{id}_U$
	give
	\[
	(s_\delta\circ\sigma)|_{Z(U)}=G_{\xi^{-1}}|_{Z(U)},
	\qquad
	\tau=(s_\delta\circ G_{\xi^{-1}})|_{Z(U)}.
	\]
	Conversely, $s_\delta\circ G_{\xi^{-1}}\in\Aut(U)$, so every
	automorphism of the displayed form extends to $U$.
	Thus $\Aut_U(Z(U))$ is the product of the subgroups induced by
	sign and graph automorphisms.
	
	For $\xi\in\Gamma$, the graph automorphism $G_\xi$ preserves
	the triangular decomposition, while $\xi$ preserves the bilinear
	form and fixes $\rho$. Hence $G_\xi$ commutes with $\pi$ on $U_0$
	and with $\gamma_{-\rho}$ on $U^0$, so
	\[
	(\gamma_{-\rho}\circ\pi)(G_\xi(z))
	=G_\xi\bigl((\gamma_{-\rho}\circ\pi)(z)\bigr)
	\qquad(z\in Z(U)).
	\]
	By Proposition~2.1 and Lemma~4.8, we may therefore compare
	the induced graph and sign actions on $(U^0_{\mathrm{ev}})^W$.
	Recall the Weyl-orbit sums
	\[
	M_\lambda=\sum_{\nu\in W\lambda}K_\nu
	\qquad(\lambda\in Q\cap2P)
	\]
	from the proof of Proposition~4.9. Since $\xi W\xi^{-1}=W$,
	\[
	G_\xi(M_\lambda)
	=\sum_{\nu\in W\lambda}K_{\xi(\nu)}
	=M_{\xi(\lambda)}.
	\]
	
	Suppose that $G_\xi|_{Z(U)}=s_\delta|_{Z(U)}$.
	Since $P/Q$ is finite, choose a positive integer $m$ such that
	$mP\subseteq Q$. For every $i$, we have $2m\omega_i\in2Q$,
	and every sign character is trivial on $2Q$. Thus
	\[
	M_{2m\xi(\omega_i)}
	=G_\xi(M_{2m\omega_i})
	=s_\delta(M_{2m\omega_i})
	=M_{2m\omega_i}.
	\]
	The monomials $K_\nu$ are linearly independent, so
	$2m\xi(\omega_i)$ and $2m\omega_i$ lie in the same Weyl orbit.
	Both are dominant, and each Weyl orbit has a unique dominant
	representative. Hence
	\[
	2m\xi(\omega_i)=2m\omega_i
	\qquad(1\leq i\leq n),
	\qquad \xi=1.
	\]
	It follows that the induced graph action is faithful and that
	its image has trivial intersection with the sign-induced subgroup.
	The factorization therefore determines $\xi$ uniquely and,
	by Proposition~4.9, also determines $\delta|_{Q\cap2P}$ uniquely.
	
	For the group structure, identify $\Gamma$ with the induced graph automorphism subgroup through the homomorphism $\xi\mapsto G_\xi|_{Z(U)}$. The defining formulas on the generators give
\[
G_\xi\circ s_\delta\circ G_\xi^{-1}=s_{\delta\circ\xi^{-1}},
\]
so the subgroup induced by sign automorphisms is normal in $\Aut_U(Z(U))$. By Proposition~4.9, this subgroup is isomorphic to $\Hom((P/Q)[2],\{\pm1\})$, with $s_\delta|_{Z(U)}$ corresponding to the character
\[
\chi(\mu+Q)=\delta(2\mu)
\qquad(\mu+Q\in(P/Q)[2]).
\]
The character corresponding to $s_{\delta\circ\xi^{-1}}|_{Z(U)}$ takes $\mu+Q$ to
\[
(\delta\circ\xi^{-1})(2\mu)
=\delta(2\xi^{-1}(\mu))
=\chi\bigl(\xi^{-1}(\mu+Q)\bigr).
\]
Thus conjugation by $G_\xi$ induces the action $(\xi\cdot\chi)(x)=\chi(\xi^{-1}x)$. Together with the factorization, faithful graph action, and trivial intersection proved above, this gives the stated semidirect-product decomposition. Under this identification, the automorphism $(s_\delta\circ G_{\xi^{-1}})|_{Z(U)}$ corresponds to $(\chi,\xi^{-1})$.
	\qedhere
\end{proof}

\begin{corollary}\label{cor:center-groups-all-types}
	For each Cartan type, the group of automorphisms induced on the center
	is given by
	\[
	\Aut_U(Z(U))\cong
	\begin{cases}
		\mathbb Z/2\mathbb Z,
		& A_n,\ n=1\text{ or }n\text{ even},\\
		(\mathbb Z/2\mathbb Z)^2,
		& A_n,\ n\geq3\text{ odd},\\
		\mathbb Z/2\mathbb Z,
		& B_n,\ C_n,\ n\geq2,\\
		(\mathbb Z/2\mathbb Z)^2,
		& D_n,\ n\geq5\text{ odd},\\
		(\mathbb Z/2\mathbb Z)^2\rtimes\mathbb Z/2\mathbb Z,
		& D_n,\ n\geq6\text{ even},\\
		S_4,
		& D_4,\\
		\mathbb Z/2\mathbb Z,
		& E_6,\ E_7,\\
		1,
		& E_8,\ F_4,\ G_2.
	\end{cases}
	\]
	For even $n\geq6$ in type $D_n$, the nonidentity element of
	$\mathbb Z/2\mathbb Z$ interchanges the two standard generators of
	$(\mathbb Z/2\mathbb Z)^2$ in the displayed semidirect product.
\end{corollary}

\begin{proof}
	By Theorem~4.10,
	\[
	\Aut_U(Z(U))
	\cong\Hom\bigl((P/Q)[2],\{\pm1\}\bigr)\rtimes\Gamma,
	\qquad
	(\xi\cdot\chi)(x)=\chi(\xi^{-1}x).
	\]
	Except in type $D_n$ with $n$ even, Proposition~4.7 shows that
	$(P/Q)[2]$ is trivial or cyclic of order two. Hence
	\[
	\Hom\bigl((P/Q)[2],\{\pm1\}\bigr)\cong
	\begin{cases}
		\{1\}, & (P/Q)[2]=0,\\
		\mathbb Z/2\mathbb Z,
		& (P/Q)[2]\cong\mathbb Z/2\mathbb Z.
	\end{cases}
	\]
	In either case, the character group has no nonidentity automorphisms.
	Thus the action of $\Gamma$ is trivial, and
	\[
	\Aut_U(Z(U))
	\cong\Hom\bigl((P/Q)[2],\{\pm1\}\bigr)\times\Gamma.
	\]
	The corresponding cases now follow from Proposition~4.7 and the
	Dynkin-diagram automorphism groups in Section~2.3.
	
	For type $D_n$ with even $n\geq6$, Bourbaki numbering gives
	\[
	(P/Q)[2]=P/Q
	=(\mathbb Z/2\mathbb Z)(\omega_{n-1}+Q)
	\oplus(\mathbb Z/2\mathbb Z)(\omega_n+Q).
	\]
	The nonidentity element $\xi\in\Gamma$ satisfies
	\[
	\xi(\omega_{n-1})=\omega_n,
	\qquad
	\xi(\omega_n)=\omega_{n-1}.
	\]
	Since $\xi^{-1}=\xi$,
	\[
	\begin{aligned}
		(\xi\cdot\chi)(\omega_{n-1}+Q)
		&=\chi(\omega_n+Q),\\
		(\xi\cdot\chi)(\omega_n+Q)
		&=\chi(\omega_{n-1}+Q).
	\end{aligned}
	\]
	Under the identification
	\[
	\chi\longmapsto
	\bigl(\chi(\omega_{n-1}+Q),\chi(\omega_n+Q)\bigr)
	\in\{\pm1\}^2,
	\]
	the action of $\xi$ therefore interchanges the two coordinates. Thus
	\[
	\Aut_U(Z(U))
	\cong(\mathbb Z/2\mathbb Z)^2\rtimes\mathbb Z/2\mathbb Z,
	\]
	with the stated action.
	
	For type $D_4$,
	\[
	(P/Q)[2]=P/Q
	=\{0,\omega_1+Q,\omega_3+Q,\omega_4+Q\},
	\]
	and $\Gamma\cong S_3$ permutes the three nonzero elements faithfully.
	Suppose that $\xi\in\Gamma$ acts trivially on the character group.
	Then
	\[
	\chi(\xi^{-1}x)=\chi(x)
	\]
	for every character $\chi$ and every $x\in(P/Q)[2]$.
	Since characters separate elements of $(P/Q)[2]$, we have
	$\xi^{-1}x=x$ for every $x\in(P/Q)[2]$.
	The action of $\Gamma$ on $(P/Q)[2]$ is faithful, so $\xi=1$.
	Thus the induced action on
	$\Hom((P/Q)[2],\{\pm1\})\cong(\mathbb Z/2\mathbb Z)^2$
	is faithful and permutes its three nonidentity elements as the full
	symmetric group. Consequently,
	\[
	\Aut_U(Z(U))
	\cong(\mathbb Z/2\mathbb Z)^2\rtimes S_3
	\cong S_4.
	\qedhere
	\]
\end{proof}

We conclude by characterizing the toral restrictions of automorphisms that
fix the center pointwise and determining the quotient
$\ker\rho_Z/\ker\rho_0$.

\begin{proposition}
	Let $\tau\in\Aut(U^0)$. Then $\tau\in\rho_0(\ker\rho_Z)$
	if and only if there exist unique $\varepsilon\in\mathcal E$
	and $w\in W$ such that
	\[
	\varepsilon|_{Q\cap2P}=1,
	\qquad
	\tau(K_\lambda)=\varepsilon(\lambda)K_{w(\lambda)}
	\quad(\lambda\in Q).
	\]
	Consequently,
	\[
	\ker\rho_Z/\ker\rho_0
	\cong
	\Hom\bigl(Q/(Q\cap2P),\{\pm1\}\bigr)\rtimes W.
	\]
	Here $\Hom\bigl(Q/(Q\cap2P),\{\pm1\}\bigr)$ is identified with $\{\delta\in\mathcal E:\delta|_{Q\cap2P}=1\}$, and the action of $W$ is given by
	\[
	(w\cdot\delta)(\lambda)=\delta(w^{-1}\lambda)
	\qquad(w\in W,\ \lambda\in Q).
	\]
\end{proposition}

\begin{proof}
	Let $\sigma\in\Aut(U)$. By Theorem~4.1, write
	\[
	\sigma(K_\lambda)=\varepsilon(\lambda)K_{\varphi(\lambda)}
	\qquad(\lambda\in Q),
	\]
	where $\varepsilon\in\mathcal E$ and $\varphi\in\Aut(\Phi)$.
	Put $\delta=\varepsilon\circ\varphi^{-1}$. By Lemma~2.4,
	$s_\delta\circ\sigma$ is sign-free on the Cartan torus, with
	\[
	(s_\delta\circ\sigma)(K_\lambda)=K_{\varphi(\lambda)}
	\qquad(\lambda\in Q).
	\]
	By Theorem~3.6, there exist unique $\xi\in\Gamma$ and $v\in W$
	such that
	\[
	D^{-1}A_{s_\delta\circ\sigma}D=[\xi]\vartheta(v).
	\]
	Proposition~4.5 and $s_\delta^2=\mathrm{id}_U$ give
	\[
	(s_\delta\circ\sigma)|_{Z(U)}
	=G_{\xi^{-1}}|_{Z(U)},
	\qquad
	\sigma|_{Z(U)}
	=(s_\delta\circ G_{\xi^{-1}})|_{Z(U)}.
	\]
	
	The argument in the proof of Theorem~4.1 gives
	\[
	(\varphi(\lambda),\mu)
	=(\lambda,(\xi v)(\mu))
	=((\xi v)^{-1}(\lambda),\mu)
	\qquad(\lambda\in Q,\ \mu\in P),
	\]
	where the second equality follows from the invariance of the
	bilinear form under $\xi v$. Since $P$ spans the rational
	weight space and the form is nondegenerate,
	\[
	\varphi=(\xi v)^{-1}.
	\]
	By the uniqueness assertion in Theorem~4.10,
	$\sigma\in\ker\rho_Z$ if and only if
	$\xi=1$ and $\delta|_{Q\cap2P}=1$.
	Since $\varphi=(\xi v)^{-1}$ and $W\cap\Gamma=\{1\}$,
	the condition $\xi=1$ is equivalent to $\varphi\in W$.
	When $\varphi\in W$, it preserves $Q\cap2P$, so
	$\delta|_{Q\cap2P}=1$ if and only if
	$\varepsilon|_{Q\cap2P}=1$. Thus
	\[
	\sigma\in\ker\rho_Z
	\quad\text{if and only if}\quad
	\varphi\in W
	\ \text{and}\ 
	\varepsilon|_{Q\cap2P}=1.
	\]
	
	Now let $\tau\in\Aut(U^0)$. If
	$\tau\in\rho_0(\ker\rho_Z)$, choose $\sigma\in\ker\rho_Z$
	with $\rho_0(\sigma)=\tau$. The criterion above gives the
	required form of $\tau$. Conversely, if $\tau$ has the
	stated form, Theorem~4.1 provides an extension
	$\sigma\in\Aut(U)$, and the same criterion gives
	$\sigma\in\ker\rho_Z$. Hence
	$\tau\in\rho_0(\ker\rho_Z)$.
	The uniqueness of $\varepsilon$ and $w$ follows from
	Theorem~4.1.
	
	By Corollary~4.6, $\ker\rho_0\subseteq\ker\rho_Z$.
	Therefore the restriction
	\[
	\rho_0|_{\ker\rho_Z}:
	\ker\rho_Z\longrightarrow\rho_0(\ker\rho_Z)
	\]
	is a surjective group homomorphism with kernel $\ker\rho_0$.
	The first isomorphism theorem gives
	\[
	\ker\rho_Z/\ker\rho_0
	\xrightarrow{\sim}\rho_0(\ker\rho_Z),
	\qquad
	\sigma\ker\rho_0\longmapsto\sigma|_{U^0}.
	\]
	
	Under the isomorphism in Proposition~2.7,
	\[
	\tau_{\varepsilon,w}
	\longmapsto(\varepsilon\circ w^{-1},w).
	\]
	Since $W$ preserves $Q\cap2P$, the condition
	$\varepsilon|_{Q\cap2P}=1$ is equivalent to
	$(\varepsilon\circ w^{-1})|_{Q\cap2P}=1$.
	The characterization proved above therefore yields
	\[
	\rho_0(\ker\rho_Z)
	\cong
	\{\delta\in\mathcal E:\delta|_{Q\cap2P}=1\}\rtimes W,
	\]
	with the action
	\[
	(w\cdot\delta)(\lambda)=\delta(w^{-1}\lambda)
	\qquad(w\in W,\ \lambda\in Q).
	\]
	Finally, a character of $Q$ is trivial on $Q\cap2P$
	if and only if it factors through $Q/(Q\cap2P)$.
	This identifies the character subgroup with
	$\Hom\bigl(Q/(Q\cap2P),\{\pm1\}\bigr)$ and gives the
	asserted isomorphism. Explicitly, under this identification,
	it is given by
	\[
	\sigma\ker\rho_0
	\longmapsto(\varepsilon\circ w^{-1},w),
	\qquad
	\sigma(K_\lambda)=\varepsilon(\lambda)K_{w(\lambda)}
	\quad(\lambda\in Q).
\qedhere	\]
\end{proof}

\section*{Acknowledgements}
This work was supported by the National Natural Science Foundation
of China (No.~12371041). The author acknowledges support from the Jiangsu Provincial Department of Education for his visiting scholarship at Nanjing University, during which this paper was written, and sincerely thanks Professor Gongxiang Liu for his warm hospitality and generous assistance throughout the visit. The author also thanks Professor Libin Li of Yangzhou University for his helpful suggestions and constructive comments.
\section*{Statements and Declarations}

The author declares that there are no competing interests relevant to the content of this article.

%
%


\end{document}